\documentclass[aap,MSNbibl,citesort,dvips]{arximspdf}
\usepackage{mathbh}
\usepackage{graphicx}

\doi{10.1214/11-AAP807} 
\volume{22}
\issue{6}
\pubyear{2012}
\firstpage{2240}
\lastpage{2281}

\makeatletter
\newcommand{\R}{\mathbb{R}}
\newcommand{\N}{\mathbb{N}}
\newcommand{\Z}{\mathbb{Z}}
\newcommand{\C}{\mathbb{C}}
\newcommand{\E}{\mathbb{E}}
\newcommand{\Prob}{\mathbb{P}}
\newcommand{\Var}{\operatorname{Var}}
\newcommand{\ind}{\mathbh{1}}
\newcommand{\eps}{\varepsilon}

\newtheorem{proposition}{Proposition}
\newtheorem{lemma}{Lemma}
\newtheorem{corollary}{Corollary}
\newtheorem{theorem}{Theorem}
\newproclaim{hyp}{Assumption}
\newproclaim{notation}{Notation}
\newproclaim{example}{Example}
\newcommand{\eqref}[1]{(\ref{#1})}
\makeatother

\begin{document}
\begin{frontmatter}

\title{Crossings of smooth shot noise processes\thanksref{T1}}
\thankstext{T1}{Supported by the French ANR Grant ``MATAIM'' No.
ANR-09-BLAN-0029-01.}
\runtitle{Crossings of smooth shot noise processes}

\begin{aug}
\author[A]{\fnms{Hermine} \snm{Bierm\'{e}}\corref{}\ead[label=e1]{hermine.bierme@mi.parisdescartes.fr}}
\and
\author[B]{\fnms{Agn\`es} \snm{Desolneux}\ead[label=e2]{agnes.desolneux@mi.parisdescartes.fr}}

\runauthor{H. Bierm\'{e} and A. Desolneux}
\affiliation{Universit\'e Paris Descartes and Universit\'{e} Fran\c
{c}ois Rabelais de Tours, and~Universit\'e~Paris~Descartes}
\address[A]{MAP5 (UMR CNRS 8145)\\
Sorbonne Paris Cit\'{e}\\
Universit\'e Paris Descartes\\
45 rue des Saints-P\`{e}res, 75006 Paris\\
France\\
and\\
LMPT (UMR CNRS 6083)\\
F\'{e}d\'{e}ration Denis Poisson\\
Universit\'{e} Fran\c{c}ois Rabelais de Tours\\
Parc de Grandmont, 37200 Tours\\
France\\
\printead{e1}} 
\address[B]{MAP5 (UMR CNRS 8145)\\
Sorbonne Paris Cit\'{e}\\
Universit\'e Paris Descartes\\
45 rue des Saints-P\`{e}res, 75006 Paris\\
France\\
\printead{e2}}
\end{aug}

\received{\smonth{5} \syear{2011}}
\revised{\smonth{9} \syear{2011}}

%
\begin{abstract}
In this paper, we consider smooth shot noise processes and their
expected number of level crossings. When the kernel response function
is sufficiently smooth,
the mean number of crossings function is obtained through an integral
formula. Moreover, as the intensity increases, or equivalently,
as the number of shots becomes larger, a normal convergence to the
classical Rice's formula for Gaussian processes is obtained.
The Gaussian kernel function, that corresponds to many applications in
physics, is studied in detail and two different regimes are exhibited.
\end{abstract}

%


\end{frontmatter}
%

\section{Introduction}
In this paper, we will consider a \textit{shot noise process} which is a
real-valued random process given by
%
\begin{equation}
X(t) = \sum_i \beta_i g(t-\tau_i) , \qquad  t\in\R,
\label{SN1deq}
\end{equation}
where $g$ is a given (deterministic) measurable function (it will be
called the \textit{kernel function} of the shot noise process), the $\{
\tau_i\}$ are the points of a Poisson point process on the line of
intensity $\lambda\nu(ds)$, where $\lambda>0$ and $\nu$ is a
positive $\sigma$-finite measure on~$\R$ and the $\{\beta_i\}$ are
independent copies of a random variable $\beta$ (called the
\textit{impulse}), independent of $\{\tau_i\}$.

Shot noise processes are related to many problems in physics as they
result from the superposition of ``shot effects'' which occur at random.
Fundamental results were obtained by Rice~\cite{Rice44}. Daley \cite
{Daley71} gave sufficient conditions on the kernel function to ensure
the convergence of the formal series in a preliminary work. General
results, including sample paths properties, were given by Rosi{\'n}ski
\cite{Rosinski} in a more general setting. In most of the literature
the measure $\nu$ is the Lebesgue measure on $\R$ such that the shot
noise process
is a stationary one. In order to derive more precise sample paths
properties and especially crossings rates, mainly two properties have
been extensively exhibited and used.
The first one is the Markov property, which is valid, choosing a
noncontinuous positive causal kernel function, that is, $0$ for negative
time. This is the case, in particular, of the exponential kernel
$g(t)=e^{-t}\ind_{t\ge0}$ for which explicit distributions and
crossings rates can be obtained~\cite{Orsingher}. A simple formula
for the expected numbers of level crossings is valid for more general
kernels of this type but resulting shot noise processes are
nondifferentiable~\cite{BarDavid,Hsing}. The infinitely divisible
property is the second main tool. Actually,
this allows us to establish convergence to a Gaussian process as the
intensity increases~\cite{Papoulis71,HeinrichSchmidt85}. Sample paths
properties of Gaussian processes have been extensively studied and fine
results are known concerning the level crossings of smooth Gaussian
processes (see~\cite{Azais,Cramer}, e.g.).

The goal of the paper is to study the crossings of a shot noise process
in the general case when the kernel function $g$ is smooth.
In this setting we lose Markov's property but the shot noise process
inherits smoothness properties. Integral formulas for the number of
level crossings of smooth processes was generalized to the non-Gaussian
case by~\cite{Leadbetter66} but it uses assumptions that rely on
properties of some densities, which may not be valid for shot noise
processes. We derive integral formulas for the mean number of crossings function
and pay a special interest in the continuity of this function with
respect to the level.
Exploiting further on normal convergence, we exhibit a Gaussian regime
for the mean number of crossings function when the intensity goes to infinity.
A particular example, which is studied in detail, concerns the shot
noise process
where $\beta=1$ almost surely and $g$ is a Gaussian kernel of width
$\sigma$,
\[
g(t)=g_{\sigma} (t) = \frac{1}{\sigma\sqrt{2\pi}} e^{-t^2/2\sigma
^2} .
\]
Such a model has many applications because it is solution of the heat
equation (we consider $\sigma$ as a variable), and it thus models a
diffusion from random sources (the points of the Poisson point
process).

The paper is organized as follows. In Section~\ref{Crossingssec}, we
consider crossings for general smooth processes. We give an explicit
formula for the Fourier transform of the mean number of crossings
function of a process $X$ in terms of the characteristic function of
$(X(t),X'(t))$.
One of the difficulties is then to obtain results for the mean number
of crossings of a given level $\alpha$ and not only for almost every
$\alpha$.
Thus we focus on the continuity property of the mean number of
crossings function.
Section~\ref{crossSN} is devoted to crossings for a smooth shot noise
process $X$ defined by \eqref{SN1deq}. In order to get the continuity
of the mean number of crossings function, we study the question of the
existence and the boundedness of a probability density for $X(t)$.
In Section~\ref{HighIntensitysec}, we show how, and in which sense,
the mean number of crossings function converges to the one of a
Gaussian process when the intensity $\lambda$ goes to infinity. We
give rates of this convergence.
Finally, in Section~\ref{Gaussiankernel}, we study in detail the case
of a Gaussian kernel of width $\sigma$. We are mainly interested in
the mean number of local extrema of this process, as a function of
$\sigma$. Thanks to the heat equation, and also to scaling properties
between $\sigma$ and $\lambda$, we prove that the mean number of
local extrema is a decreasing function of $\sigma$, and give its
asymptotics as $\sigma$ is small or large.

\section{Crossings of smooth processes}
\label{Crossingssec}
The goal of this section is to investigate crossings of general smooth
processes in order to get results for smooth shot noise processes. This
is a very different situation from the one studied in [\cite*{BarDavid}, \cite*{Orsingher,Hsing}] where shot noise processes are nondifferentiable.
However, crossings of smooth processes have been extensively studied,
especially in the Gaussian processes realm (see~\cite{Azais}, e.g.) which are second order processes. Therefore, in the whole
section, we will consider second order processes which are both almost
surely and mean square continuously differentiable (see~\cite{Adler},
Section 2.2, e.g.). This implies, in particular, that the
derivatives are also second order processes.
Moreover, most of known results on crossings are based on assumptions
on density probabilities, which are not well adapted for shot noise
processes. In this section, we revisit these results with a more
adapted point of view based on characteristic functions.

When $X$ is an almost surely continuously differentiable process on $\R
$, we can consider its multiplicity function on an interval $[a,b]$
defined by
%
\begin{equation}\label{multfunc}
\forall\alpha\in\R\qquad  N_{X}(\alpha,[a,b])=\#\{t\in
[a,b] ;   X(t)=\alpha\}.
\end{equation}
This defines a positive random process taking integer values.
Let us briefly recall some points of ``vocabulary.'' For a given level
$\alpha\in\R$, a point $t\in[a,b]$ such that $X(t)=\alpha$
is called ``crossing'' of the level $\alpha$. Then $N_X(\alpha
,[a,b])$ counts the number of crossings of the level $\alpha$ in the
interval $[a,b]$. Now we have to distinguish three different types of
crossings (see, e.g.,~\cite{Cramer}): the up-crossings that
are points for which $X(t)=\alpha$ and $X'(t)>0$, the down-crossings
that are points for which $X(t)=\alpha$ and $X'(t)<0$ and the
tangencies that are points for which $X(t)=\alpha$ and $X'(t)=0$.

Let us also
recall that according to Rolle's theorem, whatever the level $\alpha$
is,
\[
N_X(\alpha,[a,b])\le N_{X'}(0,[a,b])+1 \qquad \mbox{a.s.}
\]
Note that when there are no tangencies of $X'$ for the level $0$, then
$N_{X'}(0,[a,b])$ is the number of local extrema for $X$, which
corresponds to the sum of the number of local minima (up zero-crossings
of $X'$) and of the number of local maxima (down zero-crossings of
$X'$).

Dealing with random processes, one may be more interested in the mean
number of crossings.
We will denote by $C_X(\alpha,[a,b])$ the mean number of crossings of
the level $\alpha$ by the process $X$ in $[a,b]$,
%
\begin{equation}\quad
C_X(\alpha,[a,b]) =\E(N_{X}(\alpha,[a,b]))= \E( \#\{
t\in[a,b] \mbox{ such that } X(t)=\alpha\} ).
\end{equation}
Let us emphasize that this function is no more with integer values and
can be continuous with respect to $\alpha$. When, moreover, $X$ is a
stationary process, by the additivity of means, we get $C_X (\alpha,
[a,b])=(b-a)C_X (\alpha, [0,1])$ for $\alpha\in\R$. In this case
$C_X (\alpha, [0,1])$ corresponds to the mean number of crossings of
the level $\alpha$ per unit length. Let us also recall that when $X$
is a strictly stationary ergodic process, the ergodic theorem states
that $(2T)^{-1}N_X(\alpha,[-T,T]){\longrightarrow}_{T\rightarrow+\infty
}C_X (\alpha, [0,1])$ a.s. (see~\cite{Cramer}, e.g.).

\subsection{A Fourier approach for the mean number of crossings function}
One way to obtain results on crossings for almost every level $\alpha$
is to use the well-known \textit{co-area formula} which is, in fact, valid
in the
more general framework of bounded variations functions (see, e.g.,~\cite{EvansGariepy}). When $X$ is an almost surely
continuously differentiable process on $[a,b]$, for any bounded and
continuous function $h$ on $\R$, we have
%
\begin{equation}\label{coarea}
\int_{\R} h(\alpha) N_{X}(\alpha, [a,b])   \,d\alpha=\int_a^b
h(X(t)) |X'(t)|   \,dt \qquad \mbox{a.s.}
\end{equation}
In particular, when $h=1$, this shows that $\alpha\mapsto N_{X}(\alpha
, [a,b])$ is integrable on $\R$ and $\int_{\R} N_{X}(\alpha, [a,b])
  \,d\alpha=\int_a^b |X'(t)|   \,dt$ is the total variation of $X$ on
$[a,b]$. Moreover, taking the expected values we get by Fubini's
theorem that
\[
\int_{\R} C_X (\alpha, [a,b])   \,d\alpha=\int_a^b\E(|X'(t)|)\,dt.
\]
Therefore, when the total variation of $X$ on $[a,b]$ has finite
expectation, the function $\alpha\mapsto C_X(\alpha,[a,b])$ is
integrable on $\R$. This is the case when $X$ is also mean square
continuously differentiable since then the function $t\mapsto\E
(|X'(t)|)$ is continuous on $[a,b]$. Let us emphasize that this
implies, in particular, that $ C_X (\alpha, [a,b])<+\infty$ for
almost every level $\alpha\in\R$ but
one cannot conclude for a fixed given level. However, it allows us to
use Fubini's theorem such that, taking expectation in \eqref{coarea},
for any bounded continuous function $h$,
%
\begin{equation}
\int_{\R} h(\alpha) C_X (\alpha,[a,b])   \,d\alpha=\int_a^b\E(
h(X(t)) |X'(t)| )\,dt.
\label{SNcoareaeq}
\end{equation}
In the following theorem we obtain a closed formula for
the Fourier transform of the mean number of crossings function, which
only involves characteristic functions of the process. This can be helpful
when considering shot noise processes whose characteristic functions
are well known.

\begin{theorem} \label{Crossingspropgenerale} Let $a, b\in\R$ with
$a< b$. Let $X$ be an almost surely and mean square continuously
differentiable process on $[a,b]$. Then $\alpha\mapsto C_X(\alpha,\break
[a,b]) \in L^1(\R)$ and its Fourier transform $u\mapsto\widehat
{C_X}(u,[a,b])$ is given by
%
\begin{equation}
\widehat{C_X}(u,[a,b])=\int_a^b\E\bigl( e^{iuX(t)} |X'(t)|
\bigr)\,dt .
\label{FourierCrossingseq}
\end{equation}
Moreover, if $\psi_t$ denotes the joint characteristic function of
$(X(t),X'(t))$, then $\widehat{C_X}(u,[a,b])$ can be computed by
\begin{eqnarray*}
\widehat{C_X}(u,[a,b])
& = & - \frac{1}{\pi}\int_a^b \int_{0}^{+\infty} \frac{1}{v}
\biggl( \frac{\partial\psi_t}{\partial v} (u,v) - \frac{\partial\psi
_t}{\partial v} (u,-v) \biggr)   \,dv  \,dt\\
& = & - \frac{1}{\pi} \int_a^b\int_{0}^{+\infty} \frac
{1}{v^2}\bigl(\psi_t(u,v) + \psi_t(u,-v) - 2\psi_t(u,0) \bigr)
\,dv \, dt.
\end{eqnarray*}
\end{theorem}

\begin{pf} Choosing in equation (\ref{SNcoareaeq}) $h$ of the form
$h(x)=\exp(iux)$ for any $u$ real, shows that
$\widehat{C_X}(u,[a,b])=\int_a^b\E( e^{iuX(t)} |X'(t)|
)\,dt$. Let us now identify the right-hand term.
Let $\mu_t(dx,dy)$ denote the law of $(X(t),X'(t))$. Then the joint
characteristic function $\psi_t(u,v)$ of $(X(t),X'(t))$ is
\[
\psi_t(u,v)=\E\bigl( \exp\bigl(iuX(t)+ivX'(t)\bigr)\bigr) = \int_{\R^2}
e^{iux+ivy} \mu_t(dx,dy).
\]
Since the random vector $(X(t),X'(t))$ has moments of order two, then
$\psi_t$ is twice continuously differentiable on $\R^2$.
Now, let us consider the integral
\begin{eqnarray*}
I_A & = & \int_{0}^{A} \frac{1}{v}\biggl( \frac{\partial\psi
_t}{\partial v} (u,v) - \frac{\partial\psi_t}{\partial v} (u,-v)
\biggr)   \,dv\\
& =& \int_{v=0}^{A} \int_{x,y\in\R^2} \frac{iy
e^{iux+ivy}- iy e^{iux-ivy}}{v}\mu_t(dx,dy)  \,dv \\
& = & - 2 \int_{v=0}^{A} \int_{\R^2} y e^{iux} \frac{\sin
(vy)}{v}\mu_t(dx,dy)   \,dv\\
& =& -2 \int_{\R^2} y e^{iux} \int
_{v=0}^{Ay} \frac{\sin(v)}{v}   \,dv \mu_t(dx,dy).
\end{eqnarray*}
The order of integration has been reversed thanks to Fubini's theorem
[$|y e^{iux} \frac{\sin(vy)}{v}|\leq y^2 $ which is integrable on
$[0,A]\times\R^2$ with respect to $dv\times\mu_t(dx,dy)$, since
$X'(t)$ is a second order random variable].
As $A$ goes to $+\infty$, then $\int_{v=0}^{Ay} \frac{\sin(v)}{v}
  \,dv$ goes to $\frac{\pi}{2}\operatorname{sign}(y)$, and moreover, for
all $A$, $x$ and $y$, we have $|y e^{iux} \int_{v=0}^{Ay} \frac{\sin
(v)}{v}   \,dv| \leq3 |y| $, thus by Lebesgue's dominated convergence theorem,
the limit of $-\frac{1}{\pi}I_A$ exists as $A$ goes to infinity and
its value is
\begin{eqnarray*}
&&\lim_{A\to+\infty} -\frac{1}{\pi} \int_{0}^{A} \frac
{1}{v}\biggl( \frac{\partial\psi_t}{\partial v} (u,v) - \frac
{\partial\psi_t}{\partial v} (u,-v) \biggr)   \,dv \\
&&\qquad= \int_{\R^2}
|y| e^{iux} \mu_t(dx,dy) = \E\bigl( e^{iuX(t)} |X'(t)| \bigr).
\end{eqnarray*}
The second expression in the proposition is simply obtained by
integration by parts in the above formula.
\end{pf}

The last expression considerably simplifies when $X$ is a stationary
Gaussian process almost surely and mean square continuously
differentiable on $\R$. By independence of $X(t)$ and $X'(t)$
we get $\psi_t(u,v)=\phi_{X}(u)\phi_{X'}(v)$ where $\phi_X$,
respectively, $\phi_{X'}$, denotes the characteristic function of
$X(t)$, respectively, $X'(t)$ (independent of $t$ by stationarity).
Then, the Fourier transform of the mean number of crossings function is
given by
\[
\widehat{C_X}(u,[a,b])=- \frac{b-a}{\pi}\phi_X(u) \int_{\R} \frac
{1}{v}\frac{\partial\phi_{X'}}{\partial v} (v)   \,dv.
\]
By the inverse Fourier transform we get a weak Rice's formula
%
\begin{equation}\label{RiceGaussien}
C_X(\alpha,[a,b])=\frac{b-a}{\pi}\biggl(\frac{m_2}{m_0}
\biggr)^{1/2}e^{-(\alpha-\mathbb{E}(X(0)))^2/2m_0}\qquad  \mbox{for a.e. }\alpha\in\R,
\end{equation}
where $m_0=\operatorname{Var}(X(t))$ and $m_2=\operatorname{Var}(X'(t))$. Let us
quote that in fact Rice's formula holds for all level $\alpha\in\R$
and as soon as $X$ is a.s. continuous (see~\cite{Azais}, Exercise~3.2)
in the sense that $C_X(\alpha,[a,b])=+\infty$ if $m_2=+\infty$.

However, in general, the knowledge of $\widehat{C_X}(u,[a,b])$ only
allows us to get almost everywhere results on $C_X(\alpha,[a,b])$
itself, which can still be used in practice as explained in~\cite{Rychlik00}.

\subsection{Mean number of crossings for a given level} One way to
derive results on $C_X(\alpha,[a,b])$ for a given level $\alpha$ is
to use Kac's counting formula (see~\cite{Azais}, Lemma~3.1), which we
recall now. When $X$ is almost surely continuously differentiable on
$[a,b]$ such that for $\alpha\in\R$
%
\begin{eqnarray}\label{assKac}
\Prob\bigl(\exists t\in[a,b] \mbox{ s.t. } X(t)=\alpha\mbox{ and }
X'(t)=0\bigr)&=&0  \quad \mbox{and}
\nonumber
\\[-8pt]
\\[-8pt]
\nonumber
\Prob
\bigl(X(a)=\alpha\bigr)=\Prob\bigl(X(b)=\alpha\bigr) &=&0,
\end{eqnarray}
then,
%
\begin{equation}\label{Kac}
N_X(\alpha,[a,b])=\lim_{\delta\rightarrow0}\frac{1}{2\delta}\int
_a^b\ind_{|X(t)-\alpha|<\delta}|X'(t)|\,dt   \qquad    \mbox{a.s.}
\end{equation}
The first part of assumption \eqref{assKac} means that the number of
tangencies for the level~$\alpha$ is $0$ almost surely. The following
proposition gives a simple criterion to check this.
%
\begin{proposition}\label{NoTangencyprop} Let $a,b \in\R$ with $a\le b$. Let $X$ be a real
valued random process almost surely ${\mathcal C}^2$ on $
[a,b]$. Let us assume that there exists
$\phi\in L^1(\R)$ and $c>0$ such that
\[
\forall t\in[a,b] \qquad   \bigl|\E\bigl(e^{iu X(t)}
\bigr)\bigr|\le c\phi(u) .
\]
Then,
\[
\forall\alpha\in\R\qquad   \mathbb{P}\bigl( \exists t\in
[a,b], X(t)=\alpha\mbox{ and } X'(t)=0\bigr)=0.
\]
\end{proposition}

\begin{pf}
Let $M>0$ and let denote $A_M$ the event corresponding to
\[
\max_{t\in[a,b]} |X''(t)|\le2 M
\]
such that
$ \mathbb{P}( \exists t\in[a,b], X(t)=\alpha,
X'(t)=0)=\lim_{M\rightarrow+\infty}\mathbb{P}(\exists
t\in[a,b], X(t)=\alpha, X'(t)=0, A_M).$
Let us assume that there exists $t\in[a,b]$ such that $X(t)=\alpha$
and $X'(t)=0$. Then for any $n\in\N$ there exists
$k_n\in[2^na,2^nb]\cap\Z$ such that $|t-2^{-n}k_n|\le2^{-n}$ and,
by the first order Taylor formula,
%
\begin{equation}\label{Taylor1}
|X(2^{-n}k_n)-\alpha|\le2^{-2n}M.
\end{equation}
Therefore, let us denote
\[
B_n=\bigcup_{k_n\in[2^na,2^nb]\cap\Z}\{
|X(2^{-n}k_n)-\alpha|\le2^{-2n}M\}.
\]
Since $(B_n\cap A_M)_{n\in\N}$ is a decreasing sequence
we get
\[
\mathbb{P}\bigl(\exists t\in[a,b]; X(t)=\alpha, X'(t)=0, A_M
\bigr)\le\lim_{n\rightarrow+\infty}\mathbb{P}(B_n\cap A_M).
\]
But, according to assumption, for any $n\in\N$ the random variable
$X(2^{-n}k_n)$ admits a uniformly bounded density function.
Therefore, there exists $c'>0$ such that
\[
\mathbb{P}\bigl(|X(2^{-n}k_n)-\alpha|\le2^{-2n}M\bigr)\le c'2^{-2n}M.
\]
Hence,
$\mathbb{P}(B_n\cap A_M)\le(b-a+1)c'2^{-n}M,$
which yields the result.
\end{pf}

Now taking expectation in \eqref{Kac} gives an upper bound on
$C_X(\alpha,[a,b])$, according to Fatou's lemma,
\[
C_X(\alpha,[a,b])\le\liminf_{\delta\rightarrow0} \frac
{1}{2\delta}\int_a^b\E\bigl(\ind_{|X(t)-\alpha|<\delta
}|X'(t)|\bigr)\,dt.
\]
This upper bound is not very tractable without assumptions on the
existence of a bounded joint density for the law of $(X(t),X'(t))$.
As far as shot noise processes are concerned, one can exploit the
infinite divisibility property by considering the mean number of
crossings function
of the sum of independent processes. The next proposition gives an
upper bound in this setting. Another application of this proposition
will be seen in Section~\ref{Gaussiankernel} where we will decompose
a shot noise process into the sum of two independent processes (for
which crossings are easy to compute) by partitioning the set of points
of the Poisson process.

\begin{proposition}[(Crossings of a sum of independent processes)]\label{CrossingSumprop} Let
$a, b\in\R$ with $a< b$.
Let $n\ge2$ and $X_j$ be independent real-valued processes almost
surely and mean square two times continuously differentiable on $[a,b]$
for $1\le j\le n$. Assume that there exist constants $c_j$ and
probability measures $d\mu_j$ on $\R$ such that if $d{P}_{X_j(t)}$
denotes the law of $X_j(t)$, then
\[
\forall t\in[a,b] \qquad  d{P}_{X_j(t)}\leq c_j d\mu_j\qquad
\mbox{for } 1\le j\le n.
\]
Let $X$ be the process obtained by $  X=\sum_{j=1}^nX_j$
and assume that $X$ satisfies \eqref{assKac} for
$\alpha\in\R$. Then
%
\begin{equation}
C_X(\alpha,[a,b]) \leq\sum_{j=1}^n\biggl(\prod_{i\neq j}c_i\biggr)
\bigl(C_{X'_j}(0,[a,b]) +1\bigr).
\label{CrossingSumeq}
\end{equation}
Moreover, in the case where all the $X_j$ are stationary on $\R$,
\[
C_X(\alpha,[a,b]) \leq\sum_{j=1}^nC_{X'_j}(0,[a,b]).
\]
\end{proposition}

\begin{pf}
We first need an elementary result. Let $f$ be a $C^1$ function on
$[a,b]$, then for all $\delta>0$, and for all $x\in\R$, we have
%
\begin{equation}\label{maj0}
\frac{1}{2\delta}\int_a^b \ind_{|f(t)-x|\leq\delta} |f'(t)|   \,dt
\leq N_{f'}(0,[a,b])+1 .
\end{equation}
This result (that can be found as an exercise at the end of Chapter 3
of~\cite{Azais}) can be proved this way: let $a_1 <\cdots< a_n$
denote the points at which $f'(t)=0$ in $[a,b]$. On each interval
$[a,a_1]$, $[a_1,a_2], \ldots, [a_n,b]$, $f$ is monotonic and thus
$\int_{a_i}^{a_{i+1}}\ind_{|f(t)-x|\leq\delta} |f'(t)|   \,dt \leq
2\delta$. Summing up these integrals, we have the announced result.

For the process $X$, since it satisfies the conditions of Kac's formula
\eqref{assKac}, by \eqref{Kac} and Fatou's lemma,
\[
C_X(\alpha,[a,b]) \leq\liminf_{\delta\rightarrow0}
\frac{1}{2\delta} \int_a^b \E\bigl(\ind_{|X(t)-\alpha|\leq\delta}
|X'(t)|\bigr)   \,dt .
\]
Now, for each $\delta>0$, we have
\[
\E\bigl(\ind
_{|X(t)-\alpha|\leq\delta} |X'(t)|\bigr)\leq \sum_{j=1}^n\E\bigl(\ind
_{|X_1(t)+\cdots+X_n(t)-\alpha|\leq\delta} |X'_j(t)|\bigr).
\]
Then,
thanks to the independence of $X_1,\ldots,X_n$ and to the bound on the
laws of $X_j(t)$, we get
\begin{eqnarray*}
&&\hspace*{-4pt}\int_a^b \E\bigl(\ind_{|X_1(t)+\cdots+X_n(t)-\alpha|\leq\delta}
|X'_1(t)|\bigr)   \,dt\\
&&\hspace*{-18pt}\qquad =  \int_a^b \int_{\R^{n-1}}\E\bigl(\ind
_{|X_1(t)+x_2+\cdots+x_n-\alpha|\leq\delta} |X'_1(t)| | \\
&&\hspace*{54pt}\qquad{}
X_2(t)=x_2,\ldots,X_n(t)=x_n\bigr)  \, dP_{X_2(t)}(x_2),\ldots,
dP_{X_n(t)}(x_n)   \,dt \\
&&\hspace*{-18pt}\qquad \leq \Biggl(\prod_{j=2}^n c_j\Biggr)\int_{\R^{n-1}} \int_a^b
\E\bigl(\ind_{|X_1(t)+x_2+\cdots+x_n-\alpha|\leq\delta} |X'_1(t)|\bigr)
\,dt  \, d\mu_2(x_2),\ldots, d\mu_n(x_n).
\end{eqnarray*}
Now, \eqref{maj0} holds almost surely for $X_1$, taking expectation we get
\[
\frac{1}{2\delta}\int_a^b \E\bigl(\ind_{|X_1(t)+x_2+\cdots+x_n-\alpha
|\leq\delta} |X'_1(t)|\bigr)   \,dt \le C_{X'_1}(0,[a,b])+1.
\]
Using the fact the $d\mu_j$ are probability measures we get
\[
\frac{1}{2\delta}\int_a^b \E\bigl(\ind_{|X_1(t)+\cdots+X_n(t)-\alpha
|\leq\delta} |X'_1(t)|\bigr)   \,dt\le\biggl(\prod_{j=2}^n c_j\biggr)
\bigl(C_{X'_1}(0,[a,b])+1\bigr).
\]
We obtain similar bounds for the other terms. Since this holds for all
$\delta>0$, we have the bound \eqref{CrossingSumeq} on the expected
number of crossings of the level $\alpha$ by the process $X$.

When the $X_j$ are stationary, things become simpler; we can take
$c_j=1$ for all $1\le j\le n$, and also by stationarity we have that
for all $p\geq1$ integer
$C_X(\alpha,[a,b+p(b-a)]) = (p+1)C_X(\alpha,[a,b])$. Now, using
\eqref{CrossingSumeq} for all $p$, then dividing by \mbox{$(p+1)$}, we have
that for all $p$,
$  C_X(\alpha,[a,b]) \leq\sum_{j=1}^nC_{X'_j}(0,[a,b])
+\frac{n}{p+1} .$ Finally, letting $p$ go to infinity, we have the result.
\end{pf}

As previously seen, taking the expectation in Kac's formula only allows
us to get an upper bound for $C_X$. However, under stronger assumptions
(see~\cite{Leadbetter66}, Theorem 2), one can justify the interversion
of the limit and the expectation.
In particular, one has to assume that $(X(t),X'(t))$ admits a density
$p_t$ continuous in a neighborhood of $\{\alpha\}\times\R$.
Rice's formula states that
%
\begin{equation}\label{Ricefor}
C_X(\alpha,[a,b])=\int_a^b\int_{\R}|z|p_t(\alpha,z)   \,dz  \,dt
<+\infty,
\end{equation}
such that, under appropriate assumptions, one can prove that the mean
number of crossings function $\alpha\mapsto C_X(\alpha,[a,b])$ is
continuous on $\R$.

\section{Crossings of smooth shot noise processes} \label{crossSN}
From now on, we focus on a shot noise process $X$ given by the formal
sum \eqref{SN1deq}, which can also be written as the stochastic integral
%
\begin{equation}\label{SN1deq2}
X(t)=\int_{\R\times\R} zg(t-s)N(ds, dz),
\end{equation}
where $N$ is a Poisson random measure of intensity $\lambda\nu(ds)
F(dz)$, with $F$ the law of the impulse $\beta$ (see \cite
{Kallenberg}, Chapter 10, e.g.).
We focus in this paper on stationary shot noise processes for which
$\nu(ds)=ds$ is the Lebesgue measure. Such processes are obtained as
the almost sure limit of truncated shot noise processes defined for
$\nu_T(ds)=\ind_{[-T,T]}(s)\,ds$, as $T$ tends to infinity.
Therefore, from now on and in all the paper, the measure $\nu(ds)$ is
the Lebesgue measure $ds$ or the measure $\nu_T(ds)$.
Then, assuming that the random impulse $\beta$ is an integrable random
variable of $L^1(\Omega)$ and that the kernel function $g$ is an
integrable function of $L^1(\R)$, it is enough to ensure the almost
sure convergence of the infinite sum (see also Campbell's theorem and~\cite{HeinrichSchmidt85}). When, moreover, $\beta\in L^2(\Omega)$
and $g\in L^2(\R)$, the process $X$ defines a second order process.

\subsection{Regularity and Fourier transform of the mean number of
crossings function} Under further regularity assumptions on the kernel
function we obtain the following sample paths regularity for the shot
noise process itself.

\begin{proposition} \label{diff}Let $\beta\in L^2(\Omega)$. Let $g
\in{\mathcal C}^2(\R)$ such that $g, g', g'' \in L^1(\R)$. Then $X$
is almost surely
and mean square continuously differentiable on $\R$
with
\[
X'(t)=\sum_i \beta_i g'(t-\tau_i)\qquad  \forall t\in\R.
\]
\end{proposition}

\begin{pf}
Let $A>0$ and remark that for any $s\in\R$ and $|t|\le A$, since $g
\in{\mathcal C}^1(\R)$,
\[
|g(t-s)|=\biggl|\int_0^tg'(u-s)\,du +g(-s)\biggr|
\le\int_{-A}^A|g'(u-s)|\,du +|g(-s)|,
\]
such that by Fubini's theorem, since $g, g' \in L^1(\R)$,
\[
\int_{\R}\sup_{t\in[-A,A]}|g(t-s)| \,ds\le2A\int_{\R}|g'(s)|
\,ds+\int_{\R}|g(s)| \,ds<+\infty.
\]
Therefore, since $\beta\in L^1(\Omega)$, the series $
\sum_i \beta_i\sup_{t\in[-A,A]}|g(t-\tau_i)|$ converges
almost surely which means that $ \sum_i \beta_i g(\cdot
-\tau_i)$ converges uniformly on $[-A,A]$ almost surely. This implies
that the sample paths of $X$ are almost surely continuous on $\R$.
Similarly, since $g' \in{\mathcal C}^1(\R)$ and $g', g'' \in L^1(\R
)$, almost surely the series $ \sum_i \beta_ig'(\cdot
-\tau_i)$ converges uniformly on $[-A,A]$ and therefore $X$ is
continuously differentiable on $[-A,A]$
with $  X'(t)= \sum_i \beta_i g'(t-\tau_i)$ for all
$t\in[-A,A]$. Note that the same holds true on $[-A+n, A+n]$ for any
$n\in\Z$, which concludes for the almost sure continuous
differentiability on $\R=\bigcup_{n\in\Z}[-A+n, A+n]$.

Now, let us be concerned with the mean square continuous differentiability.
First, $g,g'\in L^1(\R)$ implies that $g\in L^\infty(\R)\cap L^1(\R
)\subset L^2(\R)$ such that $X$ is a second order process since $\beta
\in L^2(\Omega)$.
Its covariance function is given by $S(t,t')=\operatorname{Cov}
(X(t),X(t'))=\lambda\mathbb{E}(\beta^2)\int_{\R
}g(t-s)g(t'-s)\nu(ds)$.
Similarly,
we also have that $g'\in L^2(\R)\cap L^\infty(\R)$ and $X'$ is a
second order process.
According to~\cite{Adler}, Theorem~2.2.2, it is sufficient to remark
that assumptions on $g$ ensure
that $\frac{\partial^2 S }{\partial t\,\partial t'}$ exists and is
finite at any point $(t,t)\in\R^2$ with
$
\frac{\partial^2 S }{\partial t\,\partial t'}(t,t)=\lambda\mathbb
{E}(\beta^2)\int_{\R}g'(t-s)g'(t-s)\nu(ds).
$
Therefore, for all $t\in\R$, the limit
$\lim_{h \rightarrow0}\frac{X(t+h)-X(t)}{h}$
exists in $L^2(\Omega)$ and is equal to $X'(t)$ by unicity. Moreover,
the covariance function of $X'$ is given by $(t,t')\mapsto\lambda
\mathbb{E}(\beta^2)\int_{\R}g'(t-s)g'(t'-s)\nu(ds)$.
\end{pf}

Iterating this result one can obtain higher order smoothness
properties. In particular, it is straightforward to obtain the
following result for Gaussian kernels.

\begin{example*}[(Gaussian kernel)]\label{ex1}
\hspace*{-1pt}Let $\beta\in L^2(\Omega)$,
$g(t)=g_1(t)=\frac{1}{\sqrt{2\pi}}  \exp(-t^2/2)$ and $X$ given by
\eqref{SN1deq}. Then,
the process $X$ is almost surely and mean square smooth on $\R$.
Moreover, for any $n\in\N$,
\[
\forall t\in\R\qquad   X^{(n)}(t)=\sum_i \beta_i
g_1^{(n)}(t-\tau_i)=\sum_i \beta_i(-1)^n H_n (t-\tau_i) g_1(t-\tau
_i)   ,
\]
where $H_n$ is the Hermite polynomial of order $n$.
\end{example*}

From now on, in order to work with almost sure and mean square
continuously differentiable process, we make the following assumption:
{\renewcommand{\theequation}{\Alph{equation}}
\setcounter{equation}{0}
\begin{equation}\label{eqA}
 g \in{\mathcal C}^2(\R)\qquad \mbox{with } g, g', g'' \in L^1(\R) .
\end{equation}}
\indent Therefore, choosing $\beta\in L^2(\Omega)$, the shot noise process
$X$ satisfies the assumptions of
Theorem~\ref{Crossingspropgenerale} such that the Fourier transform
of its mean number of crossings function can be written with respect to
$\psi_t$, the joint characteristic function of $(X(t),X'(t))$, given
by (see~\cite{Kallenberg}, Lemma 10.2, e.g.)
%
\setcounter{equation}{14}
\begin{eqnarray}\label{psit}
\hspace*{20pt}\forall u, v\in\R \qquad   \psi_{t}(u,v) &=& \E\bigl(e^{iuX(t)+vX'(t)}\bigr)
\nonumber
\\[-4pt]
\\[-12pt]
\nonumber
&=&
\exp\biggl(\int_{\R\times\R} \bigl[e^{i z (u
g(t-s)+vg'(t-s))} - 1\bigr]  \lambda\nu(ds) F(dz)\biggr).
\end{eqnarray}

In order to get stronger results on the mean number of crossings
function we first have to investigate the existence of a density when
considering a shot noise process $X$, or more precisely, a shot noise
vector-valued process $(X,X')$. Then we
consider an $\R^d$-valued shot noise process given
on $\R$ by
%
\begin{equation}
Y(t) = \sum_i \beta_i h(t-\tau_i) ,
\label{SNndeq}
\end{equation}
where $h\dvtx \R\mapsto\R^d$ is a given (deterministic) measurable
vectorial function in $L^1(\R)$. In this setting one can recover $X$
given by \eqref{SN1deq} with $d=1$ and $h=g$, or recover
$(X,X')$ (if it exists) with $d=2$ and $h=(g,g')$. It will be
particularly helpful to see $Y$ as the almost
sure limit of a truncated shot noise process $Y_T$ defined for $\nu
_T(ds)=\ind_{[-T,T]}(s)\,ds$, as $T>0$ tends to infinity.
Therefore, from now on and in all the paper, we use the following
notation.

\begin{notation*}For any $T>0$, we denote by $Y_T$, respectively, $X_T,$
when \mbox{$d=1$}, the shot noise process given by \eqref{SNndeq},
respectively, \eqref{SN1deq}, obtained for $\nu_T(ds)=\ind
_{[-T,T]}(s)\,ds$. We simply denote by $Y$, respectively, $X,$ when
$d=1$, the shot noise process obtained for $\nu$ the Lebesgue measure.
\end{notation*}

\subsection{Existence of a density and continuity of the mean number
of crossings function}
Let us remark that for $d\ge1$ and $T>0$, the shot noise process $Y_T$
satisfies
%
\begin{equation}\label{SNndtrunc}
Y_T(\cdot)=\sum_{|\tau_i|\le T}\beta_ih(\cdot-\tau_i)\mathop{=}^\mathrm{f.d.d.}\sum_{i=1}^{\gamma_T} \beta_ih\bigl(\cdot-U_T^{(i)}\bigr),
\end{equation}
where
%
\begin{equation}\label{gammaT}
\gamma_T=\#\{i; \tau_i\in[-T,T]\}
\end{equation}
is a Poisson random variable of parameter $\lambda\nu_T(\R)=2\lambda
T$ and
$\{U_T^{(i)}\}$ are i.i.d. with uniform law on $[-T,T]$ independent
from $\gamma_T$ and $\{\beta_i\}$. Here and in the sequel the
convention is that $ \sum_{i=1}^{0}=0$ and, as usual,
$\stackrel{\mathrm{f.d.d.}}{=}$ stands for the equality in finite dimensional
distributions.

Moreover, for any $M> T$, one can write $Y_M$ as the sum of two
independent processes $Y_T$ and $Y_M-Y_T$
such that the existence of a density for the random vector $Y_T(t)$
implies the existence of a density for the random vector $Y_M(t)$ and therefore
for $Y(t)$. Note also that by stationarity $Y(s)$ will also admit a
density for any $s\in\R$. Such a remark can be used, for instance, to
establish an integral equation to compute or approximate the density in
some examples~\cite{Orsingher,LowenTeich90,Gubner96}. However, the
shot noise process may not have a density. For example,
when $h$ has compact support, there exists $A>0$ such that $h(s)=0$ for
$|s|>A$. Then, for any $T\ge A$, we get $X_T(0)=X_A(0)=X(0)$ such that
$\Prob(X_T(0)=0)=\Prob(X(0)=0)\geq\Prob(\gamma_A=0)>0$, which
proves that $X_T(0)$ and $X(0)$ don't have a density. Such a behavior
is extremely linked to the number of
points of the Poisson process $\{\tau_i\}$ that are thrown in the
interval of study. Therefore, by conditioning we obtain the following criterion.

\begin{proposition}\label{exidens} If there exists $m\ge1$ such that
for all $T>0$ large enough, conditionally on $\{\gamma_T=m\}$, the
random variable $Y_T(0)$ admits a density, then, conditionally on $\{
\gamma_T\ge m\}$, the random variable $Y_T(0)$ admits a density. Moreover,
$Y(0)$ admits a density.
\end{proposition}

\begin{pf} Let $T>0$ be large enough.
First, let us remark that conditionally on $\{\gamma_T=m\}$,
$  Y_T(0)\stackrel{d}{=} \sum_{i=1}^{m} \beta
_ih(U_T^{(i)})$. Next, notice that
if a random vector $V$ in $\R^d$ admits a density $f_V$ then, for
$U_T$ with uniform law on $[-T,T]$ and $\beta$ with law $F$,
independent of $V$, the random vector $W=V+\beta h(U_T)$ admits
$w\in\R^{d}\mapsto\frac{1}{2T}\int_{\R}\int
_{-T}^Tf_V(w-zh(t))\,dtF(dz)$ for density.
Therefore, by induction,
the assumption implies that $ \sum_{i=1}^{n}\beta_i
h(U_T^{(i)})$ has a density, for any $n\ge m$.
This proves that, conditionally on $\{\gamma_T\ge m\}$, the random
variable $Y_T(0)$ admits a density.

To prove that $Y(0)$ admits a density, we follow the same lines as in
\cite{Baccelli}, proof of Proposition A.2.
Let $A\subset \R^d$ be a Borel set with Lebesgue measure $0$,
since $Y_T(0)$ and $Y(0)-Y_T(0)$ are independent
\[
\mathbb{P}\bigl(Y(0)\in A\bigr)=\mathbb{P}\bigl(Y_T(0)+\bigl(Y(0)-Y_T(0)\bigr)\in A\bigr)=\int_{\R
^d}\mathbb{P}\bigl(Y_T(0)\in A-y\bigr)\mu_T(dy)
\]
with $\mu_T$ the law of $Y(0)-Y_T(0)$. But for any $y\in\R^d$,
\begin{eqnarray*}
\mathbb{P}\bigl(Y_T(0)\in A-y\bigr)&=&\mathbb{P}\Biggl(\sum_{i=1}^{\gamma_T}
\beta_i h\bigl(U_T^{(i)}\bigr)\in A-y\Biggr)\\
&=&\sum_{n=0}^{+\infty}\mathbb{P}\Biggl(\sum_{i=1}^{\gamma_T}
\beta_i h\bigl(U_T^{(i)}\bigr)\in A-y  \Big|  \gamma_T=n\Biggr)\mathbb
{P}(\gamma_T=n)\\
&=&\sum_{n=0}^{m-1}\mathbb{P}\Biggl(\sum_{i=1}^{n} \beta
_ih\bigl(U_T^{(i)}\bigr)\in A-y\Biggr)\mathbb{P}(\gamma_T=n),
\end{eqnarray*}
since $A-y$ has Lebesgue measure $0$ and $ \sum_{i=1}^{n}
\beta_i(h(U_T^{(i)}))$ has a density for any $n\ge m$.
Hence, for any $T>0$ large enough,
\[
\mathbb{P}\bigl(Y(0)\in A\bigr)\le\mathbb{P}(\gamma_T\le m-1).
\]
Letting $T\rightarrow+\infty$ we conclude that $\mathbb{P}(Y(0)\in
A)=0$ such that $Y(0)$ admits a~density.
\end{pf}

Let us emphasize that $Y_T(0)$ does not admit
a density since $\mathbb{P}( Y_T(0)=0)\ge \mathbb{P}( \gamma_T=0)>0$.
Let us also mention that Breton~\cite{Breton10} gives a similar assumption for
real-valued shot noise series in his Proposition 2.1.
In particular, his Corollary 2.1 can be adapted in our vector-valued setting.

\begin{corollary} \label{cordensity}  Let $h\dvtx\R\mapsto \R^d$ be an integrable
function and $\beta=1$ a.s. Let us define $h_d\dvtx\R^d\mapsto \R^d$ by $h_d(x)=h(x_1)+\cdots+h(x_d)$,
for $x=(x_1,\ldots,x_d)\in\R^d$. If the $h_d$ image measure of the $d$-dimensional
Lebesgue measure is absolutely continuous with respect to the $d$-dimensional Lebesgue measure then the
random vector $Y(0)$, given by \eqref{SNndeq}, admits a density.
\end{corollary}

\begin{pf}
Let $A\subset\R^d$ a Borel set with Lebesgue measure $0$ then the
assumptions ensure that
$\int_{\R^d}\ind_{h_d(x)\in A}\,d x=0.$
Therefore, for any $T>0$, using the notation of Proposition~\ref{exidens},
\[
\mathbb{P}\Biggl(\sum_{i=1}^d h\bigl(U_T^{(i)}\bigr)\in A\Biggr)=\frac
{1}{(2T)^d}\int_{[-T,T]^d}\ind_{h_d(x)\in A}\,d x=0.
\]
Hence, $ \sum_{i=1}^d h(U_T^{(i)})$ admits a density and
Proposition~\ref{exidens} gives the conclusion.
\end{pf}

\begin{example*}[(Gaussian kernel)]
Let $g(t)=\frac{1}{\sqrt{2\pi}}
\exp(-t^2/2)$, $\beta=1$ a.s. and $X$ given by \eqref{SN1deq}.
Let us consider $h=(g,g')$ and $h_2\dvtx (x_1,x_2)\in\R^2\mapsto
h(x_1)+h(x_2)$. The Jacobian of $h_2$
is
\[
J(h_2)(x_1,x_2)=\frac{1}{2\pi}(1+x_1 x_2)(x_1-x_2)\exp
\bigl(-(x_1^2+x_2^2)/2\bigr) .
\]
Hence, the $h_2$ image measure of the $2$-dimensional Lebesgue measure
is absolutely continuous with respect to the $2$-dimensional Lebesgue measure.
Then, for any $t\in\R$, the law of the random vector $(X(t),X'(t))$
is absolutely continuous with respect to the Lebesgue measure. Note
that, in particular, this implies the existence of a density for
$X(t)$. However, this density is not bounded (and therefore not
continuous) in a neighborhood of $0$ as proved in the following proposition.
\end{example*}

%
\begin{proposition} Let us assume for sake of simplicity that $\beta
=1$ a.s. and let $g$ denote the kernel function of the shot noise
process. Then:
\begin{longlist}[1.]
\item[1.] If $g$ is such that there exist $\alpha>1$ and $A>0$ such
that $\forall|s|>A$, $|g(s)|\leq e^{-|s|^{\alpha}}$, then $\exists
\eps_0>0$ such that $\forall0<\eps<\eps_0$,
\[
\Prob( |X(t)| \le\eps) \geq\tfrac{1}{2} e^{-2\lambda T_{\eps}}
\qquad\mbox{where } T_{\eps} \mbox{ is defined by } T_{\eps} = (-\log
\eps)^{1/\alpha}.
\]
\item[2.] If $g$ is such that there exists $A>0$ such that $\forall
|s|>A$, $|g(s)|\leq e^{-|s|}$ and if $\lambda<1/4$, then $\exists\eps
_0>0$ such that $\forall0<\eps<\eps_0$,
\[
\Prob\bigl( |X(t)| \le\eps\bigr) \geq\biggl( 1-\frac{\lambda}{(1-2\lambda
)^2}\biggr) e^{-2\lambda T_{\eps}}\quad \mbox{where } T_{\eps} \mbox{
is defined by } T_{\eps} = -\log\eps.
\]
\end{longlist}
This implies in both cases that $\Prob( |X(t)| \le\eps) / \eps$
goes to $+\infty$ as $\eps$ goes to $0$, and thus the density of
$X(t)$ (if it exists) is not bounded in a neighborhood of $0$.
\end{proposition}

\begin{pf}
We start with the first case.
Let $\eps>0$ and let $T_{\eps} = (-\log\eps)^{1/\alpha}$. Assume
that $\eps$ is small enough to have $T_{\eps}>A$.
We have by definition $  X(t)\stackrel{d}{=}X(0)\stackrel
{d}{=} \sum_{i} g(\tau_i)$. If we denote $  X_{T_{\eps
}}(0)=\sum_{|\tau_i|\leq T_{\eps}} g(\tau_i)$ and $
R_{T_{\eps}}(0)=\sum_{|\tau_i|> T_{\eps}} g(\tau_i)$, then
$X_{T_{\eps}}(0)$ and $R_{T_{\eps}}(0)$ are independent and
$X(0)=X_{T_{\eps}}(0)+R_{T_{\eps}}(0)$. We also have $\Prob(|X(0)|
\le\eps) \geq\Prob(|X_{T_{\eps}}(0)|=0 \mbox{ and } |R_{T_{\eps
}}(0)| \le\eps) = \Prob(|X_{T_{\eps}}(0)|=0) \times\Prob
(|R_{T_{\eps}}(0)| \le\eps)$.
Now, on the one hand, we have $\Prob(|X_{T_{\eps}}(0)|=0)
\geq\Prob~(\mbox{there are no } \tau_i \mbox{ in } [-T_{\eps},T_{\eps}]) =
e^{-2\lambda T_{\eps}}$.
On the other hand, the first moments of the random variable $R_{T_{\eps
}}(0)$ are given by $\E(R_{T_{\eps}}(0))=\break\lambda\int_{|s|>T_\eps
}^{+\infty} g(s)   \,ds $ and $\Var(R_{T_{\eps}}(0))= \lambda\int
_{|s|> T_\eps}^{+\infty} g^2(s)   \,ds$.
Now, we use the following inequality on the tail of $\int
e^{-s^{\alpha}}$:
\[
\forall T>0 \qquad  e^{-T^{\alpha}} = \int_T^{+\infty} \alpha s^{\alpha
-1}e^{-s^{\alpha}}   \,ds \geq \alpha T^{\alpha-1}\int_T^{+\infty}
e^{-s^{\alpha}}   \,ds .
\]
Thus, we obtain bounds for the tail of $\int g$ and of $\int g^2$,
\[
\int_T^{+\infty} e^{-s^{\alpha}}   \,ds \le\frac{e^{-T^{\alpha
}}}{\alpha T^{\alpha-1}} \quad\mbox{and}\quad \int_T^{+\infty} (
e^{-s^{\alpha}})^2   \,ds \le\frac{e^{-2T^{\alpha}}}{2\alpha
T^{\alpha-1}} .
\]
Back to the moments of $R_{T_{\eps}}(0)$, since $T_{\eps}=(-\log\eps
)^{1/\alpha}$ we have
\[
|\E(R_{T_{\eps}}(0))| \le\frac{2\lambda\eps}{\alpha T_{\eps
}^{\alpha-1}} \quad\mbox{and}\quad \Var(R_{T_{\eps}}(0)) \le\frac{\lambda
\eps^2}{\alpha T_{\eps}^{\alpha-1}} .
\]
We can take $\eps$ small enough in such a way that we can assume
that\break
$|\E(R_{T_{\eps}}(0))|<\eps$. Then, using Chebyshev's inequality, we have
\begin{eqnarray*}
\Prob\bigl(|R_{T_{\eps}}(0)| \le\eps\bigr) & = & \Prob\bigl(-\eps-\E(R_{T_{\eps
}}(0)) \le R_{T_{\eps}}(0) -\E(R_{T_{\eps}}(0)) \le\eps-\E
(R_{T_{\eps}}(0))\bigr) \\
&\ge& 1-\Prob\bigl(|R_{T_{\eps}}(0) - \E(R_{T_{\eps}}(0))| \ge\eps-
|\E(R_{T_{\eps}}(0))|\bigr) \\
& \ge& 1 - \frac{ \Var(R_{T_{\eps}}(0))}{(\eps- |\E(R_{T_{\eps
}}(0))|)^2} \ge1- \frac{\lambda}{\alpha T_{\eps}^{\alpha-1}(1 -
2\lambda/\alpha T_{\eps}^{\alpha-1})^2}
\end{eqnarray*}
which is larger than $1/2$ for $T_{\eps}$ large enough (i.e., for
$\eps$ small enough).

For the second case, we can make exactly the same computations by
setting $\alpha=1$ and get $\Prob(|R_{T_{\eps}}(0)| \le\eps) \ge
1- \lambda/(1-2\lambda)^2$, which is $>0$ when $\lambda<1/4$.
\end{pf}

Such a feature is particularly bothersome when considering crossings of
these processes since most of known results
are based on the existence of a bounded density for each marginal of
the process. However, this is again linked to the number of
points of the Poisson process $\{\tau_i\}$ that are thrown in the
interval of study. By conditioning, the characteristic functions are
proved to be integrable such that conditional laws have continuous
bounded densities. The main tool is Proposition~\ref{PhaseStatiolem}
(postponed to the \hyperref[Appendix]{Appendix}) established using the classical stationary
phase estimate for oscillatory integrals (see, e.g.,~\cite{Stein}).

\begin{proposition}\label{BoundDensityShotprop} Let us assume for sake of simplicity that $\beta
=1$ a.s.,
let \mbox{$T>0$}, $a<b$ and assume that $g\in L^1(\R)$ is a function of class
$\mathcal{C}^2$ on $[-T+a,T+b]$ such that
%
\begin{eqnarray}\label{ps1}
  m&=&\min_{s\in[-T+a,T+b]}\sqrt{g'(s)^2+g''(s)^2} >0
\quad\mbox{and }
\nonumber
\\[-8pt]
\\[-8pt]
\nonumber
n_0&=&\#\{s\in[-T+a,T+b] \mbox{ s.t. } g''(s)=0 \}
<+\infty.
\end{eqnarray}
Then, conditionally on $\{\gamma_T\ge k_0\}$ with $k_0\ge3$, for all
$t\in[a,b]$ and $M\ge T$, the law of
$X_M(t)$ admits a continuous bounded density. Therefore, for any $t\in
\R$, the law of $X(t)$, conditionally on $\{\gamma_T\ge k_0\}$,
admits a continuous bounded density.
\end{proposition}

\begin{pf} Actually, we will prove that conditionally on $\{\gamma
_T\ge k_0\}$, the law of the truncated process
$  X_T(t)=\sum_{|\tau_i|\le T}g(t-\tau_i)$ admits a
continuous bounded density for $t\in[a,b]$. The result will follow,
using the fact that for $M\ge T$,
$X_M(t)=X_T(t)+(X_M(t)-X_T(t))$, with $X_M(t)-X_T(t)$ independent of
$X_T(t)$. So let us denote $\psi_{t,k_0}^T$ the characteristic
function of $X_T(t)$ conditionally on $\{\gamma_T\ge k_0\}$. Then, for
all $u\in\R$, we get
\begin{eqnarray*}
\psi_{t,k_0}^T(u)&=&\frac{1}{\mathbb{P}(\gamma_T\ge k_0)} \sum
_{k\ge k_0}\mathbb{E}\bigl(e^{iuX_T(t)}| \gamma_T=k\bigr)\mathbb
{P}(\gamma_T=k) \\
&=&\frac{1}{\mathbb{P}(\gamma_T\ge k_0)} \sum_{k\ge k_0}
\biggl(\frac{1}{2T}\int_{-T}^Te^{iug(t-s)}\,ds\biggr)^ke^{-2\lambda T}\frac
{(2\lambda T)^k}{k!}.
\end{eqnarray*}
Therefore,
%
\begin{equation}\label{psiTk0}
|\psi_{t,k_0}^T(u)|\le(2T)^{-k_0} \biggl|\int
_{-T+t}^{T+t}e^{iug(s)}\,ds \biggr|^{k_0}.
\end{equation}
Hence, using Proposition~\ref{PhaseStatiolem} on $[-T+t,T+t]\subset
[-T+a,T+b]$, one can find $C$ a positive constant that depends on $T$,
$k_0$, $\lambda$, $m$ and $n_0$ such that
for any $|u|>1/m$
\[
|\psi_{t,k_0}^T(u)|\le C |u|^{-k_0/2}.
\]
Then, since $k_0\ge3$, $\psi_{t,k_0}^T$ is integrable on $\R$ and
thanks to Fourier inverse theorem it is the characteristic function of
a bounded continuous density.~%
\end{pf}

Using similar ideas we obtain the following result concerning the
continuity of the mean number of crossings function.

\begin{theorem}\label{contcrosssn} Assume for sake of simplicity
that $\beta=1$ a.s. and that $g$ is a function of class ${\mathcal
C}^4$ on $\R$ satisfying (\ref{eqA}). Let $T>0$, $a\le b$ and assume that
for all $s\in[-T+a,T+b]$, the matrice $\Phi(s)=\bigl(
{g'(s) \atop g''(s)} \enskip{ g''(s) \atop g^{(3)}(s)}\bigr)$ and its component-wise derivative $\Phi'(s)=\bigl(
{g''(s) \atop g^{(3)}(s) }\enskip{ g^{(3)}(s) \atop g^{(4)}(s)}\bigr)$ are invertible.
Then, conditionally on $\{\gamma_T\ge k_0\}$ with $k_0\ge8$, for all
$M\ge T$, the mean number of crossings function
$\alpha\mapsto\mathbb{E}(N_{X_M}(\alpha,[a,b])| \gamma_T\ge
k_0)$
is continuous on $\R$. Moreover,
\[
\mathbb{E}\bigl(N_{X_M}(\alpha,[a,b])| \gamma_T\ge k_0
\bigr)\mathop{\longrightarrow}_{M\rightarrow+\infty}\mathbb{E}
\bigl(N_{X}(\alpha,[a,b])| \gamma_T\ge k_0\bigr)
\]
uniformly on $\alpha\in\R$.
\end{theorem}

\begin{pf}
The result follows from Rice's formula. To establish it we use~\cite
{Leadbetter66}, Theorem 2, and thus we have to check assumptions (i) to (iii) related to joint densities.
Let $t\in[a,b]$ and $M\ge T$. We write
$X_M(t)=X_T(t)+(X_M(t)-X_T(t))$ with $X_M-X_T$ independent of $X_T$. We
adopt the convention that $X_\infty=X$.
Let us write for $M\in[T,+\infty]$ and $\varepsilon$ small enough
\[
\psi_{t,\varepsilon,k_0}^M=\psi_{t,\varepsilon,k_0}^T\psi
_{t,\varepsilon}^{T,M}
\]
with $\psi_{t,\varepsilon,k_0}^M$ the characteristic function of
$(X_M(t), (X_M(t+\varepsilon)-X_M(t))/\varepsilon)$, conditionally on
$\{\gamma_T\ge k_0\}$.
Note that, $X_M-X_T$ is independent of $\gamma_T$ such that $\psi
_{t,\varepsilon}^{T,M}$ is just the characteristic function of
$(X_M(t)-X_T(t), ((X_M-X_T)(t+\varepsilon)-(X_M-X_T)(t))/\varepsilon
)$. First we prove that there exists $C>0$ such that, for all $0\le
j\le3$, for all $M\ge T$ and $\varepsilon>0$ small enough,
%
\begin{equation}\label{majopsi}
\biggl|\frac{\partial^j}{\partial v^j}\psi_{t,\varepsilon
,k_0}^M(u,v)\biggr|\le C\bigl(1+\sqrt{u^2+v^2}\bigr)^{-(k_0-3)/2}.
\end{equation}
Let us remark that, since $g', g'' \in L^1(\R)$ by (\ref{eqA}), one has
$g, g'\in L^{\infty}(\R)$. It implies, in particular, that
$g,g'\in L^1(\R)\cap L^2(\R)\cap L^3(\R)$ such that the above
partial derivatives exist.
Moreover, by Leibniz formula, for $0\le j\le3$,
one has
%
\begin{equation}\label{Leibnitz}
\frac{\partial^j}{\partial v^j}\psi_{t,\varepsilon,k_0}^M(u,v)=\sum
_{l=0}^j \pmatrix{j\vspace*{2pt}\cr l}\frac{\partial^l}{\partial v^l}\psi_{t,\varepsilon
,k_0}^T(u,v)\frac{\partial^{j-l}}{\partial v^{j-l}}\psi
_{t,\varepsilon}^{T,M}(u,v).
\end{equation}
On the one hand,
\[
\biggl|\frac{\partial^{j-l}}{\partial v^{j-l}}\psi_{t,\varepsilon
}^{T,M}(u,v)\biggr|\le\mathbb{E}\biggl(\biggl|\frac
{(X_M-X_T)(t+\varepsilon)-(X_M-X_T)(t)}{\varepsilon}
\biggr|^{j-l}\biggr)
\]
with
\[
\biggl|\frac{(X_M-X_T)(t+\varepsilon)-(X_M-X_T)(t)}{\varepsilon
}\biggr|\le\sum_{T<|\tau_i|\le M}|g_\varepsilon(t-\tau_i)|,
\]
where $g_\varepsilon(s)=\frac{1}{\varepsilon}\int_0^\varepsilon g'(s+x)\,dx$
is such that
$g_\varepsilon\in L^{\infty}(\R)\cap L^1(\R)$ with $\|g_\varepsilon
\|_{\infty}\le\|g'\|_{\infty}$ and $\|g_\varepsilon\|_{1}\le\|g'\|
_{1}$. Then, using the moment formula established in~\cite{Bassan},
one can find
$c>0$ such that for all $0\le j\le3$, with $(j-1)_+=\max(0,j-1)$,
%
\begin{eqnarray}\label{majoreste}
\biggl|\frac{\partial^{j}}{\partial v^{j}}\psi_{t,\varepsilon
}^{T,M}(u,v)\biggr|&\le&\E\biggl(\biggl(\sum_{T<|\tau_i|\le
M}|g_\varepsilon(t-\tau_i)|\biggr)^{j}\biggr)
\nonumber
\\[-8pt]
\\[-8pt]
\nonumber
&\le& c \max(1,\|g'\|
_{\infty})^{(j-1)_+}\max(1,\lambda\|g'\|_1)^j.
\end{eqnarray}
On the other hand,
\begin{eqnarray*}
\mathbb{P}(\gamma_T\ge k_0)\psi_{t,\varepsilon,k_0}^T(u,v)&=&\sum
_{k\ge k_0}\mathbb{E}\bigl(e^{iuX_T(t)+iv(X_T(t+\varepsilon
)-X_T(t))/\varepsilon}| \gamma_T=k\bigr)\mathbb{P}(\gamma_T=k)\\
&=&\sum_{k\ge k_0}\chi_{t,\varepsilon}^T(u,v)^k\mathbb{P}(\gamma_T=k),
\end{eqnarray*}
where
\[
\chi_{t,\varepsilon}^T(u,v)=(2T)^{-1}\int
_{-T+t}^{T+t}e^{iug(s)+ivg_\varepsilon(s)}\,ds
\]
is the characteristic function of $(g(t-U_T),g_\varepsilon(t-U_T))$,
with $U_T$ a uniform random variable on $[-T,T]$. It follows that
$|\chi_{t,\varepsilon}^T(u,v)|\le1$,
so that one can find $c>0$ such that for all $0\le j\le3$,
\begin{eqnarray*}
\biggl|\frac{\partial^{j}}{\partial v^{j}}\psi_{t,\varepsilon
,k_0}^T(u,v)\biggr|
&\le& c \max(1,\|g'\|_{\infty})^{(j-1)_+}\max
(1,\lambda\|g'\|_1)^j\\
&&{}\times \frac{\mathbb{P}(\gamma_T\ge k_0-j)}{\mathbb
{P}(\gamma_T\ge k_0)}|\chi_{t,\varepsilon}^T(u,v)|^{k_0-j}.
\end{eqnarray*}
This, together with \eqref{majoreste} and \eqref{Leibnitz}, implies
that one can find
$c>0$ such that for all $0\le j\le3$,
%
\begin{eqnarray}\label{majopsicond}
\biggl|\frac{\partial^{j}}{\partial v^{j}}\psi_{t,\varepsilon
,k_0}(u,v)\biggr|&\le& c \max(1,\|g'\|_{\infty})^{(j-1)_+}\max
(1,\lambda\|g'\|_1)^j
\nonumber
\\[-8pt]
\\[-8pt]
\nonumber
&&{}\times \frac{\mathbb{P}(\gamma_T\ge k_0-j)}{\mathbb
{P}(\gamma_T\ge k_0)}|\chi_{t,\varepsilon}^T(u,v)|^{k_0-j}.
\end{eqnarray}
Moreover,
let $\Phi_\varepsilon(s)=\bigl(
{g'(s) \atop g''(s)}\enskip {g_\varepsilon'(s) \atop g_\varepsilon
''(s)}%
\bigr)$ and $\Phi'_\varepsilon(s)=\bigl(
{ g''(s) \atop g^{(3)}(s)}\enskip{ g_\varepsilon''(s) \atop
g_\varepsilon^{(3)}(s)}
\bigr)$.
Then $\operatorname{det} \Phi_\varepsilon(s)$ converges to $\operatorname{det} \Phi
(s)$ as $\varepsilon\rightarrow0$, uniformly in $s\in[-T-a,T+b]$.
The assumption on $\Phi$ ensures that one can find $\varepsilon_0$
such that for
$\varepsilon\le\varepsilon_0$, the matrix $\Phi_\varepsilon(s)$ is
invertible for all $s\in[-T-a,T+b]$. The same holds true for $\Phi
'_\varepsilon(s)$.
Denote $  m=\min_{s\in[-T-a,T+b], \varepsilon\le
\varepsilon_0}\Vert \Phi_\varepsilon(s)^{-1}\Vert ^{-1}>0$,
where $\|\cdot\|$ is the matricial norm induced by the Euclidean one.
According to Proposition~\ref{PhaseStatiolem} with $n_0=0$,
\begin{eqnarray*}
&&\forall(u,v)\in\R^2 \mbox{ s.t. } \sqrt{u^2+v^2}>\frac{1}{m} ,\\
&&\qquad
 |\chi_{t,\varepsilon}^T(u,v)
|=(2T)^{-1}\biggl| \int_{-T+t}^{T+t}e^{iug(s)+ivg_\varepsilon(s)}\,ds
\biggr| \leq\frac{24 \sqrt{2} }{\sqrt{m\sqrt{u^2+v^2}}}.
\end{eqnarray*}
Therefore, one can find a constant $c_{k_0}>0$ such that,
for all $0\le j\le3$, $|\frac{\partial^{j}}{\partial
v^{j}}\psi_{t,\varepsilon,k_0}^{T}(u,v)|$ is less than
\begin{eqnarray*}
&&c_{k_0}(2T)^{-k_0+3} \max(1,\|g'\|_{\infty})^{(j-1)_+}\max(1,\lambda
\|g'\|_1)^j\\
&&\qquad{}\times \frac{\mathbb{P}(\gamma_T\ge k_0-j)}{\mathbb{P}(\gamma
_T\ge k_0)}\bigl(1+\sqrt{u^2+v^2}\bigr)^{-(k_0-3)/2}.
\end{eqnarray*}
Letting $\varepsilon$ tend to $0$ we obtain the same bounds as \eqref
{majopsi} for $\psi_{t,k_0}^{M}$ the characteristic function of
$(X_M(t),X_M'(t))$
conditionally on $\{\gamma_T\ge k_0\}$. Since $k_0\ge8$, \eqref
{majopsi} for $j=0$ ensures that $\psi_{t,\varepsilon,k_0}^M\in
L^1(\R^2)$, respectively, $\psi_{t,k_0}^M\in L^1(\R^2)$, such that,
conditionally on $\{\gamma_T\ge k_0\}$, $(X_M(t), (X_M(t+\varepsilon
)-X_M(t))/\varepsilon)$, respectively, $(X_M(t),X_M'(t))$, admits
$p_{t,\varepsilon,k_0}^M(x,z)=\frac{1}{4\pi^2}\int_{\R
^2}e^{-ixu-izv}\psi_{t,\varepsilon,k_0}^M(u,\break v)\,du \, dv$, respectively,
$p_{t,k_0}^M=\frac{1}{4\pi^2}\int_{\R^2}e^{-ixu-izv}\psi
_{t,k_0}^M(u,v)\,du\, dv$, as density.\break Moreover:

\begin{longlist}[(iii)]
\item[(i)] $p_{t,\varepsilon,k_0}^M(x,z)$ is continuous in $(t,x)$ for each
$z, \varepsilon$, according to Lebesgue's dominated convergence
theorem using the fact that $X_M$ is almost surely continuous on $\R
$.

\item[(ii)] Since $X_M$ is almost surely continuously differentiable on $\R
$ we clearly have for any $(u,v)\in\R^2$, $\psi_{t,\varepsilon
,k_0}^M(u,v)\rightarrow\psi_{t,k_0}^M(u,v)$ as $\varepsilon
\rightarrow0$. Then by Lebes\-gue's dominated convergence theorem, using
\eqref{majopsi} for $j=0$ we check that\break $p_{t,\varepsilon
,k_0}^M(x, z)\rightarrow p_{t,k_0}^M(x,z)$ as $\varepsilon\rightarrow
0$, uniformly in $(t,x)$ for each $z\in\R$.\vspace*{1pt}

\item[(iii)] For any $z\neq0$, integrating by parts we get
\[
p_{t,\varepsilon,k_0}^M(x,z)=\frac{i}{4\pi^2 z^{3}}\int_{\R
^2}e^{-ixu-izv}\frac{\partial^3}{\partial v^3}\psi_{t,\varepsilon
,k_0}^M(u,v)\,du\, dv,
\]
such that by \eqref{majopsi} for $j=3$, we check that
$p_{t,\varepsilon,k_0}^M(x,z)\le C h(z)$ for all $t, \varepsilon, x$
with $h(z)=(1+|z|^3)^{-1}$ satisfying $\int_{\R}|z|h(z)\,dz<+\infty$
and $C$ a positive constant.

Therefore,~\cite{Leadbetter66}, Theorem 2, implies that
\[
\mathbb{E}\bigl(N_{X_M}(\alpha,[a,b])| \gamma_T\ge k_0\bigr)=\int
_a^b\int_{\R}|z|p_{t,k_0}^M(\alpha,z)   \,dz\,  dt,
\]
which concludes the proof, using $p_{t,k_0}^M(\alpha,z)=\frac{i}{4\pi
^2 z^{3}}\int_{\R^2}e^{-i\alpha u-izv}\frac{\partial^3}{\partial
v^3}\times\break \psi_{t,k_0}^M(u,v)\,du\, dv$ and \eqref{majopsi} for
$j=3$.\quad\qed
\end{longlist}
\noqed\end{pf}

Note that, despite that we have closed forms, these crossings formulas
are not very tractable for general shot noise processes. However,
as the intensity $\lambda$ of the shot noise process $X$ tends to
infinity, due to its infinitely divisible property and since it is of
second order, we obtain, after renormalization, a Gaussian process at
the limit. It is then natural to hope for the same kind of asymptotics
for the mean number of crossings function. This behavior is studied in
detail in the next section.

\section{High intensity and Gaussian field}
\label{HighIntensitysec}

\subsection{General feature}
It is well known that, as the intensity $\lambda$ of the Poisson
process goes to infinity, the shot noise process converges to a normal
process. Precise bounds on the distance between the law of $X(t)$ and
the normal distribution are given by Papoulis~\cite{Papoulis71}.
Moreover, Heinrich and Schmidt~\cite{HeinrichSchmidt85} give
conditions of normal convergence for a wide class of shot noise
processes (not restricted to processes defined on $\R$, nor to Poisson
processes). In this section we obtain a stronger result for smooth
stationary shot noise processes by considering
convergence in law in the space of continuous functions. In all of this
section we continue to assume that $X$ is a stationary shot noise
process obtained for $\nu$ the Lebesgue measure on $\R$,
and we will denote $X_\lambda$ the strictly stationary shot noise
process given by \eqref{SN1deq} with intensity $\lambda>0$. Let us
define the normalized shot noise process
%
\begin{equation}\label{normalized}
Z_\lambda(t)= \frac{1}{\sqrt{\lambda}}\bigl(X_\lambda(t)-\mathbb
{E}(X_\lambda(t))\bigr),\qquad t\in\R.
\end{equation}
Then, we obtain the following result.
%
\begin{proposition}\label{CLT} Let $\beta\in L^2(\Omega)$ and $g$
satisfying (\ref{eqA}). Then,
\[
Y_\lambda=\pmatrix{ Z_\lambda\vspace*{2pt}\cr Z'_\lambda
}\mathop{\longrightarrow}_{\lambda\rightarrow+\infty}^\mathrm{f.d.d.}\sqrt{\E(\beta^2)}
\pmatrix{B\vspace*{2pt}\cr B'},
\]
where $B$ is a stationary centered Gaussian process almost surely and
mean square continuously differentiable, with covariance function
\[
\operatorname{Cov}(B(t),B(t'))=\int_{\R}g(t-s)g(t'-s)\,ds.
\]
When, moreover, $g''\in L^p(\R)$ for $p>1$,
the convergence holds in distribution on the space of continuous
functions on compact sets endowed with the topology of the uniform convergence.
\end{proposition}

\begin{pf}
We begin with the proof of the finite dimensional distributions
convergence. Let $k$ be an integer with $k\ge1$ and let
$t_1,\ldots,t_k\in\R$ and $w_1=(u_1,v_1),\ldots,w_k=(u_k,v_k) \in
\R^2$.

Let us write
\[
\sum_{j=1}^k Y_\lambda(t_j)\cdot w_j =\frac{1}{\sqrt{\lambda
}}\biggl(\sum_{i} \beta_i\widetilde{g}(\tau_i)-\E\biggl(\sum_{i}
\beta_i\widetilde{g}(\tau_i)\biggr)\biggr),
\]
for $ \widetilde{g}(s)=\sum_{j=1}^k(
u_jg(t_j-s)+v_jg'(t_j-s))$.
Therefore
\[
\log\E\bigl(e^{i\sum_{j=1}^{k} Y_\lambda
(t_j)\cdot w_j} \bigr) = \lambda\int_{\R\times\R} \biggl(
e^{iz({\widetilde{g}(s)}/{\sqrt{\lambda}})} - 1 -i
z \frac{\widetilde{g}(s)}{\sqrt{\lambda}} \biggr)  \,dsF(dz).
\]
Note that as $\lambda\rightarrow+\infty$,
\[
\lambda\biggl( e^{iz({\widetilde{g}(s)}/{\sqrt{\lambda
}})} - 1 -i z \frac{\widetilde{g}(s)}{\sqrt{\lambda}}
\biggr)\rightarrow
-\frac{1}{2}z^2\widetilde{g}(s)^2,
\]
with, for all $\lambda>0$,
\[
\biggl|\lambda\exp\biggl(iz\biggl(\frac{\widetilde{g}(s)}{\sqrt
{\lambda}}\biggr) - 1 -i z \frac{\widetilde{g}(s)}{\sqrt{\lambda
}} \biggr)\biggr| \le
\frac{1}{2}z^2 \widetilde{g}(s)^2.
\]
By the dominated convergence theorem, since $\widetilde{g}\in L^2(\R
)$ and $\beta\in L^2(\Omega)$, we get that, as $\lambda\rightarrow
+\infty$,
\[
\E\Biggl(\exp\Biggl({i\sum_{j=1}^{k} Y_\lambda
(t_j)\cdot w_j} \Biggr)\Biggr)\rightarrow
\exp\biggl(-\frac{1}{2} \E(\beta^2)\int_{\R} \widetilde
{g}(s)^2\,ds\biggr).
\]
Let us identify the limiting process. Let us recall that
$X_\lambda$ is a second order process with covariance function given by
$\operatorname{Cov}(X_\lambda(t),X_\lambda(t'))=\lambda\E(\beta
^2)S(t-t')$ with
$S(t)=\int_{\R}g(t-s)g(-s)\,ds.$
Hence, one can define $B$ to be a stationary Gaussian centered process
with $(t,t')\mapsto S(t-t')$ as covariance function.
The assumptions on $g$ ensure that the function $S$ is twice
differentiable. Therefore $B$ is mean square differentiable with $B'$ a
stationary Gaussian centered process with $(t,t')\mapsto
-S''(t-t')=\int_{\R}g'(t-t'-s)g'(-s)\,ds$ as covariance function.
Moreover,
\[
\E\bigl(\bigl(B'(t)-B'(t')\bigr)^2\bigr)=2\bigl(S''(0)-S''(t-t')\bigr)\le
2\|g'\|_{\infty}\|g''\|_1|t-t'|,
\]
such that by~\cite{Adler}, Theorem 3.4.1., the process $B'$ is almost
surely continuous on $\R$. Therefore, as in~\cite{Doob}, page 536,
one can check that almost surely
$B(t)=B(0)+\int_0^tB'(s)\,ds,$ such that $B$ is almost surely
continuously differentiable.
We conclude for the f.d.d. convergence by noticing that
\[
\int_{\R} \widetilde{g}(s)^2\,ds=\operatorname{Var}\Biggl(\sum_{j=1}^k
u_jB(t_j)+v_jB'(t_j)\Biggr).
\]
Let us prove the convergence in distribution on the space of continuous
functions on compact sets endowed with the topology of the uniform convergence.
It is enough to prove the tightness of the sequence $(Y_\lambda
)_{\lambda}$
according to~\cite{Kallenberg}, Lemma~14.2 and Theorem 14.3. Let $t,s
\in\R$ and remark that for any $q\ge1$, on the one hand,
\begin{eqnarray*}
\mathbb{E}\bigl(\bigl(Z_\lambda(t)-Z_\lambda(t')\bigr)^2\bigr)
&=&\E(\beta
^2)\int_{\R}\bigl(g(t-s)-g(t'-s)\bigr)^2\,ds\\
&\le&\E(\beta^2)\|g'\|
_{q}\|g'\|_1|t-t'|^{2-1/q}.
\end{eqnarray*}
On the other hand,
\begin{eqnarray*}
\mathbb{E}\bigl(\bigl(Z'_\lambda(t)-Z'_\lambda(t')\bigr)^2\bigr)&=&\E(\beta
^2)\int_{\R}\bigl(g'(t-s)-g'(t'-s)\bigr)^2\,ds\\
&\le&\E(\beta^2)\|
g''\|_{q}\|g''\|_1|t-t'|^{2-1/q}.
\end{eqnarray*}
Note that, assuming that $g''\in L^p(\R)$, it allows us to choose
$q=p>1$ in the second upper bound such that $2-1/q>1$. Moreover,
assumption (\ref{eqA}) implies that $g'\in L^\infty(\R)\cap L^1(\R
)\subset L^p(\R)$ such that one can also choose $q=p$ in the first
upper bound. Then, $(Y_\lambda)_{\lambda}$ satisfies a
Kolmogorov--Chentsov criterion which implies its tightness according to
\cite{Kallenberg}, Corollary 14.9.
\end{pf}


In particular, when $a<b$, the functional $(f,g)\mapsto \int_a^bh(f(t))|g(t)|\,dt $ is clearly continuous and bounded on
${\mathcal C}([a,b],\R)\times {\mathcal C}([a,b],\R)$ for any  continuous bounded function $h$  on $\R$. Then,
Proposition~\ref{CLT} implies that
\[
\int_a^b\mathbb{E}(h(Z_\lambda(t))|Z'_\lambda(t)|)\,dt \mathop{\longrightarrow}_{\lambda\rightarrow+\infty}
\int_a^b\mathbb{E}(h(B(t))|B'(t)|)\,dt.
\]
By the co-area formula \eqref{coarea}, this means the weak convergence of the mean number  of crossings function, that is,
\[
C_{Z_\lambda}(\cdot,[a,b])\rightharpoonup_{\lambda \rightarrow +\infty} C_{B}(\cdot,[a,b]) .
\]
This implies also the pointwise convergence of Fourier transforms. Such a result can be compared to the classical central limit theorem.
Numerous improved results can be obtained under stronger assumptions than the classical ones. This is the case, for instance,
for the rate of convergence derived by the Berry--Esseen theorem or the convergence of the densities.
We refer to~\cite{Feller2}, Chapters 15 and 16. Adapting the technical proofs allows us to get similar results for crossings
in the next section.

\subsection{High intensity: rate of convergence  for the mean number  of crossings function}
Let us remark that only $\E(\beta^2)$ appears in the limit field.
For sake of simplicity we may assume that $\beta=1 \mbox{ a.s.}$
Note that, according to Rice's formula~\cite{Cramer},  as  recalled in equation \eqref{RiceGaussien}, since the limit Gaussian field
 is stationary, $C_B(\alpha,[a,b])=(b-a)C_B(\alpha,[0,1])$ with
\[
C_B(\alpha,[0,1])=\frac{1}{\pi}\biggl(\frac{m_2}{m_0}\biggr)^{1/2}e^{-\alpha^2/2m_0} \qquad\forall\alpha\in\R,
\]
where  $m_0=\operatorname{Var}(B(t))=\int_{\R}g(s)^2\,ds$ and
$m_2=\operatorname{Var}(B'(t))=\int_{\R}g'(s)^2\,ds$. Moreover, its
Fourier transform is given by $\widehat{C_{B}}(u,[0,1]) = \sqrt{\frac{2m_2}{\pi}} e^{-m_0 u^2/2}$.
We obtain the following rate of convergence, for which the proof is postponed to the \hyperref[Appendix]{Appendix}.

\begin{proposition}
\label{BoundsConvprop} Let   $\beta=1$ a.s. and let $g$ satisfy (\ref{eqA}).
There exist three constants $a_1$, $a_2$ and $a_3$ (depending only on
$g$ and its derivative) such that
\begin{eqnarray*}
&&\forall \lambda >0, \ \forall u\in \R \mbox{ such that } |u|<a_1\sqrt{\lambda} \mbox{ then } \\
&&\qquad\biggl| \widehat{C_{Z_{\lambda}}}(u,[0,1]) - \sqrt{\frac{2m_2}{\pi}} e^{-m_0 u^2/2}\biggr | \leq \frac{a_2 + a_3 |u|}{\sqrt{\lambda}},
\end{eqnarray*}
where $m_0=\int_{\R}g(s)^2\,ds$ and $m_2=\int_{\R}g'(s)^2\,ds$.
\end{proposition}

Let us emphasize that this implies the uniform convergence of the Fourier transform of
the mean number  of crossings functions on any fixed interval. Moreover, taking $u=0$,
the previous upper bound may be a bit refined
such that the following corollary is in force.

\begin{corollary}\label{totvar} Let   $\beta=1$ a.s. and let $g$ satisfy (\ref{eqA}).
 The mean total variation   of the process satisfies
 \[
\forall \lambda>0 \qquad \biggl| \frac{\E(|X'_{\lambda}(t)|)}{\sqrt{\lambda}} - \sqrt{\frac{2m_2}{\pi}}\biggr| \leq \frac{14 m_3}{3\pi m_2\sqrt{\lambda}}  ,
\]
where $m_2=\int_{\R}g'(s)^2\,ds$ and $m_3=\int_{\R} |g'(s)|^3 \, ds$.
\end{corollary}

Under additional assumptions we obtain the following uniform
convergence for the mean number of crossings function. The proof is
inspired by~\cite{Feller2}, Theorem 2, page 516, concerning the
central limit theorem for densities.
%
\begin{theorem}\label{convunifgausstheo} Let $\beta=1$ a.s. Let us assume, moreover, that $g$
is a function of class ${\mathcal C}^4$ on $\R$ satisfying (\ref{eqA}) such that
for all $s\in[-1,2]$, $\Phi(s)=\bigl(
{g'(s)\atop  g''(s)}\enskip{  g''(s) \atop g^{(3)}(s)}\bigr)$ and $\Phi'(s)=\bigl(
{g''(s) \atop g^{(3)}(s) }\enskip{ g^{(3)}(s) \atop g^{(4)}(s)}\bigr)$ are invertible.

Let $\gamma_{\lambda}=\#\{i; \tau_{\lambda,i}\in[-1,1]\}$ with $\{
\tau_{\lambda,i}\}_i$ the points of a Poisson point process with
intensity $\lambda>0$.

Then
\begin{eqnarray}
C_{Z_\lambda}(\alpha,[0,1]|\gamma_\lambda\ge\lambda) \mathop{\longrightarrow}_{\lambda\rightarrow+\infty}C_{B}(\alpha
,[0,1])=\frac{1}{\pi}\biggl(\frac{m_2}{m_0}\biggr)^{1/2}e^{-\alpha
^2/2m_0}\nonumber\\
   \eqntext{\mbox{uniformly in }  \alpha\in\R,}
\end{eqnarray}
where $m_0=\int_{\R}g(s)^2\,ds$ and $m_2=\int_{\R}g'(s)^2\,ds$.
\end{theorem}

\begin{pf} Let $\lambda\ge8$.
Then, according to Theorem~\ref{Crossingspropgenerale}, $\widehat
{C_{Z_\lambda}}(u,[0,1] | \gamma_\lambda\ge\lambda)$
and $\widehat{C_B}(u,[0,1])$ are integrable such that ${C_{Z_\lambda
}}(\alpha,[0,1]|\gamma_\lambda\ge\lambda)$
and ${C_B}(\alpha,[0,1])$ are bounded continuous functions with, for
any $\alpha\in\R$,
\begin{eqnarray*}
&&|{C_{Z_\lambda}}(\alpha,[0,1] | \gamma_\lambda\ge\lambda
)-{C_B}(\alpha,[0,1])|\\
&&\qquad\le\frac{1}{2\pi}\int_{\R}
|\widehat{C_{Z_\lambda}}(u,[0,1] | \gamma_\lambda\ge
\lambda)-\widehat{C_B}(u,[0,1])|\,du.
\end{eqnarray*}
Let $u\in\R$, then
\begin{eqnarray*}
&&\widehat{C_{Z_\lambda}}(u,[0,1])-\widehat{C_{Z_\lambda}}(u,[0,1]
| \gamma_\lambda\ge\lambda)\\
&&\qquad=\frac{1}{\mathbb{P}(\gamma
_\lambda\ge\lambda)}
\E\bigl(e^{iuZ_\lambda(0)}|Z'_\lambda(0)|\ind_{\gamma_\lambda
<\lambda}\bigr)-\frac{\mathbb{P}(\gamma_\lambda<\lambda
)}{\mathbb{P}(\gamma_\lambda\ge\lambda)}\widehat{C_{Z_\lambda}}(u,[0,1]).
\end{eqnarray*}
Note that $|\widehat{C_{Z_\lambda}}(u,[0,1])|\le\E
( |Z'_\lambda(0)| )$, which is bounded according to
Corollary~\ref{totvar}, while by the Cauchy--Schwarz inequality,
\[
\bigl|\E\bigl(e^{iuZ_\lambda(0)}|Z'_\lambda(0)|\ind_{\gamma
_\lambda<\lambda}\bigr)\bigr|\le\E( Z'_\lambda(0)^2
)^{1/2}\mathbb{P}(\gamma_\lambda<\lambda)^{1/2},
\]
with $ \E( Z'_\lambda(0)^2 )=\operatorname{Var}(Z'_\lambda
(0))\le\max(1,\|g'\|_\infty)\|g'\|_1$.
Therefore, one can find $c_1>0$ such that
\[
|\widehat{C_{Z_\lambda}}(u,[0,1])-\widehat{C_{Z_\lambda
}}(u,[0,1] | \gamma_\lambda\ge\lambda)|\le c_1 \frac
{\mathbb{P}(\gamma_\lambda<\lambda)^{1/2}}{\mathbb{P}(\gamma
_\lambda\ge\lambda)}.
\]
According to Markov's inequality,
\[
\mathbb{P}(\gamma_\lambda< \lambda)=\mathbb{P}\bigl(e^{-\ln
(2)\gamma_\lambda} > e^{-\ln(2)\lambda}\bigr)\le\mathbb{E}
\bigl(e^{-\ln(2)(\gamma_\lambda-\lambda)}\bigr)=\exp\bigl(-\bigl(1-\ln
(2)\bigr)\lambda\bigr).
\]
Choosing $\lambda$ large enough such that, in particular, $\frac
{\mathbb{P}(\gamma_\lambda<\lambda)^{1/2}}{\mathbb{P}(\gamma
_\lambda\ge\lambda)}\le\frac{1}{\sqrt{\lambda}}$, according to
Proposition~\ref{BoundsConvprop} one can find $c_2$ such that
for all $|u|<\lambda^{1/8}$,
\[
|\widehat{C_{Z_{\lambda}}}(u,[0,1] | \gamma_\lambda) - \widehat
{C_B}(u,[0,1])| \leq c_2\lambda^{-3/8}.
\]
Thus we may conclude that
\[
\int_{|u|< \lambda^{1/8}}
|\widehat{C_{Z_\lambda}}(u,[0,1] | \gamma_\lambda\ge
\lambda)-\widehat{C_B}(u,[0,1])|\,du \mathop{\longrightarrow}_{\lambda
\rightarrow+\infty} 0.
\]
Now, let us be concerned with the remaining integral for $|u|\ge
\lambda^{1/8}$.
According to Theorem~\ref{Crossingspropgenerale},
\[
\widehat{C_{Z_\lambda}}(u,[0,1] | \gamma_\lambda\ge\lambda
)=\frac{e^{-iu\sqrt{\lambda}\int_{\R}g}}{\sqrt{\lambda}}\widehat
{C_{X_\lambda}}\biggl(\frac{u}{\sqrt{\lambda}},[0,1] | \gamma
_\lambda\ge\lambda\biggr),
\]
with $\widehat{C_{X_\lambda}}(\frac{u}{\sqrt{\lambda
}},[0,1] | \gamma_\lambda\ge\lambda)=\int_0^1\E(
e^{i({u}/{\sqrt{\lambda}})X_\lambda(t)} |X'_\lambda(t)| |
\gamma_\lambda\ge\lambda) \,dt $ and
\begin{eqnarray*}
&&\E\bigl( e^{i({u}/{\sqrt{\lambda}})X_\lambda(t)} |X'_\lambda
(t)| | \gamma_\lambda\ge\lambda\bigr)\\
&&\qquad= - \frac{1}{\pi}\int
_{0}^{+\infty} \frac{1}{v}\biggl( \frac{\partial\psi_{t,\lambda
}}{\partial v} \biggl(\frac{u}{\sqrt{\lambda}},\frac{v}{\sqrt
{\lambda}}\biggr) - \frac{\partial\psi_{t,\lambda}}{\partial v}
\biggl(\frac{u}{\sqrt{\lambda}},-\frac{v}{\sqrt{\lambda}}\biggr)
\biggr)   \,dv,
\end{eqnarray*}
where $\psi_{t,\lambda}$ is the characteristic function of
$(X_\lambda(t),X'_\lambda(t))$ conditionally on $\{\gamma_\lambda
\ge\lambda\}$.
Integrating by parts we obtain
\begin{eqnarray*}
&&\int_{0}^{1} \frac{1}{v}\biggl( \frac{\partial\psi_{t,\lambda
}}{\partial v} \biggl(\frac{u}{\sqrt{\lambda}},\frac{v}{\sqrt
{\lambda}}\biggr) - \frac{\partial\psi_{t,\lambda}}{\partial v}
\biggl(\frac{u}{\sqrt{\lambda}},-\frac{v}{\sqrt{\lambda}}\biggr)
\biggr)   \,dv
\\
&&\qquad=
-\frac{1}{\sqrt{\lambda}}\int_{0}^{1} \ln(v)\biggl( \frac
{\partial^2 \psi_{t,\lambda}}{\partial v^2}\biggl(\frac{u}{\sqrt
{\lambda}},\frac{v}{\sqrt{\lambda}}\biggr) - \frac{\partial^2
\psi_{t,\lambda}}{\partial v^2} \biggl(\frac{u}{\sqrt{\lambda
}},-\frac{v}{\sqrt{\lambda}}\biggr) \biggr) \,dv.
\end{eqnarray*}
Then,
according to \eqref{majopsicond}, one can find a positive constant
$c_3>0$ such that
\begin{eqnarray*}
&&\bigl|\E\bigl( e^{i({u}/{\sqrt{\lambda}})X_\lambda(t)}
|X'_\lambda(t)|  | \gamma_\lambda\ge\lambda\bigr)\bigr|\\
&&\qquad\le
c_3 \lambda^2\frac{\mathbb{P}(\gamma_\lambda\ge\lambda
-2)}{\mathbb{P}(\gamma_\lambda\ge\lambda)}\\
& &\quad\qquad{}\times\int_{\R} \biggl|\chi_{t}\biggl(\frac{u}{\sqrt{\lambda
}},\frac{v}{\sqrt{\lambda}}\biggr)\biggr|^{\lambda-2}\biggl(\frac
{1}{\sqrt{\lambda}}|\ln(|v|)|\ind_{0\le|v|\le1}+|v|^{-1}\ind
_{|v|\ge1}\biggr)\,dv,
\end{eqnarray*}
where $\chi_{t}(u,v) =\frac{1}{2}\int
_{-1+t}^{1+t}e^{iug(s)+ivg'(s)}\,ds$ is the characteristic function of
$(g(t-U),g'(t-U))$, with $U$ a uniform random variable on $[-1,1]$.
Then,
\begin{eqnarray*}
&&\int_{|u|\ge\lambda^{1/8}}
|\widehat{C_{Z_\lambda}}(u,[0,1] | \gamma_\lambda\ge
\lambda)-\widehat{C_B}(u,[0,1])|\,du \\
&&\qquad\le \int_{|u|\ge
\lambda^{1/8}}|\widehat{C_{Z_\lambda}}(u,[0,1] | \gamma
_\lambda\ge\lambda)|\,du +\int_{|u|\ge\lambda^{1/8}}
|\widehat{C_B}(u,[0,1])|\,du\\
&&\qquad= I_1(\lambda)+I_2(\lambda).
\end{eqnarray*}
Now, for $\theta\in[0,2\pi]$, let us consider the random variable
$V_{t,\theta}=\cos(\theta)g(t-U)+\sin(\theta)g'(t-U)$
such that for any $r>0$, $\chi_{t}(r\cos(\theta),r\sin(\theta
))=\mathbb{E}(e^{irV_{t,\theta}}):=\varphi_{t,\theta}(r)$. By a
change of variables in polar coordinates, since $\lambda>1$, we get
\[
I_1(\lambda)\le c_{4}(\lambda)\int_{\lambda^{1/8}}^{+\infty}\int
_0^{2\pi}\biggl| \varphi_{t,\theta}\biggl(\frac{r}{\sqrt{\lambda
}}\biggr)\biggr|^{\lambda-2} r\bigl(|\ln(r|\sin(\theta
)|)|+1\bigr)\,d\theta \,dr,
\]
with $c_4(\lambda)=c_3 \lambda^{3/2}\frac{\mathbb{P}(\gamma
_\lambda\ge\lambda-2)}{\mathbb{P}(\gamma_\lambda\ge\lambda)}$.
Since $\operatorname{det}\Phi(s)\neq0$
for any $s\in[-1+t,1+t]$, we have the following property (see \cite
{Feller2}, page 516): there exists $\delta>0$ such that
\begin{eqnarray*}
|\varphi_{t,\theta}(r)|&\le& e^{-({\kappa(t)}/{4})r^2} \qquad\forall r\in(0,\delta],\ \forall\theta\in[0,2\pi]
 \quad\mbox{and }\\
 \eta&=&\sup_{r>\delta,\theta\in[0,2\pi
]}|\varphi_{t,\theta}(r)|<1,
\end{eqnarray*}
with $\kappa(t)=\min_{\theta\in[0,2\pi]}\operatorname
{Var}(V_{t,\theta})>0$. Note also that according to Proposition~\ref
{PhaseStatiolem},
$|\varphi_{t,\theta}(r)|\le24\sqrt{\frac
{2}{m}}r^{-1/2}$ for any $r>m$ with $m=\min_{s\in[-1,2]}\|
\Phi(s)^{-1}\|^{-1}$, which may be assumed\vadjust{\goodbreak} to be larger than $\delta$.
Then, for $\lambda$ large enough such that $\lambda^{1/8}\in
(e,\delta\sqrt{\lambda})$,
\begin{eqnarray*}
I_1(\lambda)&\le& c_5(\lambda)\Biggl(\int_{\lambda^{1/8}}^{\delta
\sqrt{\lambda}}e^{-({\kappa(t)}/{8})\lambda^{1/4}}r\ln(r) \,dr+
\int_{\delta\sqrt{\lambda}}^{m\sqrt{\lambda}}\eta^{\lambda
-2}r\ln(r)\,dr\\
&&\hspace*{68pt}\qquad{}+\Biggl(24\sqrt{\frac{2}{m}}\Biggr)^5\int_{m\sqrt
{\lambda}}^{+\infty}\eta^{\lambda-7}r^{-3/2}\ln(r)\,dr\Biggr)
\end{eqnarray*}
with $c_5(\lambda)=c_4(\lambda)(\int_0^{2\pi}(2+|\ln
(|\sin(\theta)|)|)\,d\theta)$. This enables us to
conclude that
$I_1(\lambda)\longrightarrow_{\lambda\rightarrow+\infty}0$. This concludes the proof since clearly $I_2(\lambda)
\longrightarrow_{\lambda\rightarrow+\infty}0$.
\end{pf}

Notice that to obtain the convergence in Theorem \ref
{convunifgausstheo} without the conditioning on $\{\gamma_{\lambda
}\ge\lambda\}$ (which is an event of probability going to $1$
exponentially fast as $\lambda$ goes to infinity), one simply needs to
have an upper-bound polynomial in $\lambda$ on the second moment of
the number of crossings $N_{Z_\lambda}(\alpha,[0,1])$.

\section{The Gaussian kernel}\label{Gaussiankernel}

In this section we will be interested in a real application of shot
noise processes in physics. Indeed, each time a physical model is given
by sources that produce each a potential in such a way that the global
potential at a point is the sum of all the individual potentials, then
this can be modeled as a shot noise process. In particular, we will be
interested here in the temperature produced by sources of heat.
Assuming that the sources are randomly placed as a Poisson point
process of intensity $\lambda$ on the real line $\R$, then the
temperature after a time ${\sigma^2}$ on the line is given by the
following shot noise process $X_{\lambda,\sigma}$:
\[
t\in\R\mapsto X_{\lambda,\sigma}(t)= \sum_i \frac{1}{\sigma\sqrt
{2\pi}} e^{-(t-\tau_i)^2/2\sigma^2} ,
\]
where the $\{\tau_i\}$ are the points of a Poisson process of intensity
$\lambda>0$ on $\R$.
In the following, we will denote by $g_{\sigma}$ the Gaussian kernel
of width $\sigma$ defined for all $t\in\R$ by
\[
g_{\sigma} (t) = \frac{1}{\sigma\sqrt{2\pi}} e^{-t^2/2\sigma^2} .
\]

We will be interested in the crossings of $X_{\lambda,\sigma}$
because they provide information on the way the temperature is
distributed on the line. The number of local extrema of $X_{\lambda
,\sigma}$ is also interesting for practical applications since it
measures the way the temperature fluctuates on the line.
In a first part, we will be interested in the crossings of $X_{\lambda
,\sigma}$ when $\lambda$ and $\sigma$ are fixed, and then, in a
second part, we will study how the number of crossings evolves when
these two parameters change. From the point of view of applications,
this amounts to describing the fluctuations of the temperature on the
line when the time (recall that $\sigma^2$ represents the time)
increases, or when the number of sources changes.

\subsection{\texorpdfstring{Crossings and local extrema of $X_{\lambda,\sigma}$}\
{Crossings and local extrema of X lambda, sigma}}

We assume in this subsection that $\lambda>0$ and $\sigma>0$ are fixed.
Since the Gaussian kernel $g_\sigma$, and its derivatives are smooth
functions which belong to all $L^{p}$ spaces, many results of the
previous sections about crossings can be applied here. In particular,
we have:
\begin{itemize}
\item the function $\alpha\mapsto C_{X_{\lambda,\sigma}}(\alpha
,[a,b])$ belongs to $L^1(\R)$ (by Theorem~\ref{Crossingspropgenerale});
\item for any $T>0$, the function $\alpha\mapsto C_{X_{\lambda,\sigma
}}(\alpha,[a,b] | \gamma_T\geq8)$ is continuous (by Theorem \ref
{contcrosssn}), with $\gamma_T=\#\{\tau_i\in[-T,T]\}$.
\end{itemize}
This second point comes from the fact that the Gaussian kernel
satisfies the hypothesis of Theorem~\ref{contcrosssn}. Indeed, the
derivatives of $g_\sigma$ are given by $g_\sigma^{(k)}(s)=\frac
{1}{\sigma\sqrt{2\pi}}e^{-s^2/2\sigma^2} \cdot\frac
{(-1)^k}{\sigma^k}H_k(\frac{s}{\sigma})$, where the $H_k$'s are the
Hermite polynomials ($H_1(x)=x$ ; $H_2(x)= x^2-1$; $H_3(x) = x^3-3x $
and $H_4(x)=x^4-6x^2+3$). Thus, using the notation of Theorem \ref
{contcrosssn}, we get
$\det\Phi(s) = \frac{-1}{\sigma^4}(\frac{s^2}{\sigma^2}+1)
( \frac{1}{\sigma\sqrt{2\pi}}e^{-s^2/2\sigma^2} )^2 <0$ and
$\det\Phi'(s) = \frac{-1}{\sigma^6}(\frac{s^4}{\sigma^4}+3)
( \frac{1}{\sigma\sqrt{2\pi}}e^{-s^2/2\sigma^2} )^2 <0$.
These two matrices are thus invertible for all $s\in\R$.

The first point implies that for almost every $\alpha\in\R$, the
expected number of crossings of the level $\alpha$ by $X_{\lambda
,\sigma}$ is finite.
We will now prove in the following proposition that in fact, for every
$\alpha\in\R$, $ C_{X_{\lambda,\sigma}}(\alpha,[a,b])<+\infty$,
by considering the zero-crossings of the derivative $X'_{\lambda
,\sigma}$ and using Rolle's theorem.

In the sequel, we will denote by $\rho(\lambda,\sigma)$ the mean
number of
local extrema of $X_{\lambda,\sigma}$ in the interval $[0,1]$.
It is the mean number of local extrema per unit length.

\begin{proposition}\label{bound19}
We have
\[
\Prob\bigl( \exists t\in[0,1] \mbox{ such that } X'_{\lambda,\sigma
}(t)=0 \mbox{
and } X''_{\lambda,\sigma}(t)=0 \bigr) = 0,
\]
which implies that the local extrema of $X_{\lambda,\sigma}$ are
exactly the points where the derivative vanishes; in other words $\rho
(\lambda,\sigma) = \E(N_{X'_{\lambda,\sigma}}(0,[0,1]))$.
Moreover, we have the following bounds:
\[
\forall\alpha\in\R \qquad  C_{X_{\lambda,\sigma}}(\alpha
,[0,1]) \leq\rho(\lambda,\sigma)\leq\bigl(3\lambda(2+2\sigma) +1\bigr)
e^{\lambda} .
\]
\end{proposition}

\begin{pf}
For the first part of the proposition, we use Proposition \ref
{PhaseStatiolem} (in the \hyperref[Appendix]{Appendix}) with the kernel function
$h=g'_{\sigma}$ on the interval $[-T+1,T]$ for $T>0$. For this
function we can compute\vspace*{-2pt} $h'(s)=\frac{1}{\sigma^3\sqrt{2\pi
}}(-1+\frac{s^2}{\sigma^2})e^{-s^2/2\sigma^2}$ and $h''(s)=\frac
{1}{\sigma^4\sqrt{2\pi}}(3\frac{s}{\sigma}-\frac{s^3}{\sigma
^3})e^{-s^2/2\sigma^2}$, and thus $n_0=3$ and $m(\sigma,T)=\break \min
_{s\in[-T,T+1]}\sqrt{h'(s)^2+h''(s)^2}>0$ (we do not need to have an
exact value for it but notice that it is of the order of
$e^{-T^2/2\sigma^2}$ when $T$ is large). Finally, as in \eqref
{psiTk0}, we get that there is a constant $c(T,\sigma)$ which depends
continuously on $\sigma$ and $T$ such that
\[
\bigl|\E\bigl(e^{iuX'_{\lambda,\sigma}(t)} | \gamma_T\ge3\bigr)\bigr| \leq \frac{
c(T,\sigma)^3}{(1+\sqrt{|u|})^3},
\]
with $\gamma_T=\#\{\tau_i\in[-T,T]\}$.
We can now use Proposition~\ref{NoTangencyprop} and we get that for
all $T>1$,
\[
\Prob\bigl( \exists t\in[0,1] \mbox{ such that } X'_{\lambda,\sigma}(t)=0
\mbox{ and } X''_{\lambda,\sigma}(t)=0 | \gamma_T\ge3\bigr) = 0.
\]
Since the events $\{\gamma_T\ge3\}$ are an increasing sequence of
events such that $\Prob(\gamma_T\ge3)$ goes to $1$ as $T$ goes to
infinity, we obtain that
$\Prob(\exists t\in[0,1] \mbox{ such that}\break X'_{\lambda,\sigma
}(t)=0 \mbox{ and } X''_{\lambda,\sigma}(t)=0 ) = 0$.

For the second part of the proposition,
the left-hand inequality is simply a consequence of Proposition \ref
{CrossingSumprop} for the process $X'_{\lambda,\sigma}$ and $n=1$.

To obtain the right-hand inequality [the bound on $\rho(\lambda
,\sigma)$], we will apply Proposition~\ref{CrossingSumprop} to the
process $X'_{\lambda,\sigma}$ for the crossings of the level $0$ on the
interval $[0,1]$. We already know by the first part of the proposition
and by Corollary~\ref{cordensity} that condition~\eqref{assKac} for
Kac's formula is satisfied by $X'_{\lambda,\sigma}$. Then we write, for
all $t\in[0,1]$,
\begin{eqnarray*}
X'_{\lambda,\sigma}(t) = \sum_{\tau_i\in\R} g'_{\sigma}(t-\tau
_i) &= &\frac{1}{\sigma\sqrt{2\pi}} \sum_{\tau_i\in[-\sigma
,1+\sigma]} \frac{-(t-\tau_i)}{\sigma^2} e^{-(t-\tau_i)^2/2\sigma
^2}\\
&&{} + \frac{1}{\sigma\sqrt{2\pi}} \sum_{\tau_i\in\R\setminus
[-\sigma,1+\sigma]} \frac{-(t-\tau_i)}{\sigma^2} e^{-(t-\tau
_i)^2/2\sigma^2}.
\end{eqnarray*}
Let $Y_1(t)$ [resp., $Y_2(t)$] denote the first (resp., second) term.
We then have
\[
Y'_2(t)= \frac{1}{\sigma\sqrt{2\pi}} \sum_{\tau_i\in\R
\setminus[-\sigma,1+\sigma]} \biggl( \frac{(t-\tau_i)^2}{\sigma
^4} -\frac{1}{\sigma^2} \biggr) e^{-(t-\tau_i)^2/2\sigma^2}.
\]
Since $(t-\tau_i)^2>\sigma^2$ for all $t\in[0,1]$ and all $\tau
_i\in\R\setminus[-\sigma,1+\sigma]$, we get $Y'_2(t)>0$ on
$[0,1]$ and thus $ N_{Y'_2}(0,[0,1])=0$ a.s. Note that when the event
$\#\{\tau_i\in[-\sigma,1+\sigma]\}=0$ holds, then $X'_{\lambda
,\sigma}=Y_2$ such that $N_{X'_{\lambda,\sigma}}(0,[0,1])\le1$. On
the other hand, let us work conditionally on $\#\{\tau_i\in[-\sigma
,1+\sigma]\}\ge1$. The probability of this event is $1-e^{-\lambda
(1+2\sigma)}$. To study the zero-crossings of $Y'_1$, we first need an
elementary lemma.

\begin{lemma}\label{Nbracineslem}
Let $n\geq1$ be an integer. Let $P_1,\ldots, P_n$ be $n$ real
nonzero polynomials and let $a_1,\ldots, a_n$ be $n$ real numbers, then
\[
\#\Biggl\{ t\in\R\mbox{ such that } \sum_{i=1}^n P_i(t) e^{a_i t}=0 \Biggr\}
\leq\sum_{i=1}^n \deg(P_i) + n - 1 .
\]
\end{lemma}

This elementary result can be proved by induction on $n$. For $n=1$, it
is obviously true. Assume the result holds for $n\geq1$, then we prove
it for $n+1$ in the following way. For $t\in\R$, $\sum_{i=1}^{n+1}
P_i(t) e^{a_i t}=0 \Longleftrightarrow f(t):= P_{n+1}(t) + \sum
_{i=1}^{n} P_i(t) e^{(a_i-a_{n+1}) t}=0$. Let $k$ denote the degree of
$P_{n+1}$. Thanks to Rolle's theorem, we have that $N_f(0,\R)\leq
N_{f'}(0,\R) +1 \leq N_{f''}(0,\R) +2 \leq\cdots\leq
N_{f^{(k+1)}}(0,\R) + k+1$. But $f^{(k+1)}$ can be written as
$f^{(k+1)}(t) = \sum_{i=1}^{n} Q_i(t) e^{(a_i-a_{n+1}) t}$, where the
$Q_i$ are polynomials of degree $\deg(Q_i)\leq\deg(P_i)$. Thus by
induction $N_{f^{(k+1)}}(0,\R)\leq\sum_{i=1}^n \deg(P_i) + n - 1$,
and then $N_f(0,\R)\leq\sum_{i=1}^n \deg(P_i) + n - 1+k+1 \leq\sum
_{i=1}^{n+1} \deg(P_i) + n$. This proves the result for $n+1$.

Thanks to this lemma, we get that $N_{Y'_1}(0,[0,1])\leq3\#\{\tau
_i\in[-\sigma,1+\sigma] \}-1$
such that
\[
\mathbb{E}\bigl(N_{Y'_1}(0,[0,1])|\#\{\tau_i\in[-\sigma,1+\sigma]\}\ge
1\bigr)\le3\lambda(1+2\sigma)/\bigl(1-e^{-\lambda(1+2\sigma)}\bigr)-1.
\]
To use Proposition~\ref{CrossingSumprop}, we need to obtain uniform
bounds on the laws of $Y_1(t)$ and of $Y_2(t)$ when $t\in[0,1]$. As in
the notation of the proposition, we will denote these constants by
$c_1$ and $c_2$. Let us start with $Y_1$.
Let $U$ be a random variable following the uniform distribution on
$[-1-\sigma,1+\sigma]$. For $t\in[0,1]$, we can write $U$ as $U=\eta
_t U_t + (1-\eta_t) V_t$, where $U_t$ is uniform on $[-1-\sigma
+t,\sigma+t]$, $V_t$ is uniform on $[-1-\sigma,-1-\sigma+t]\cup
[\sigma+t,\sigma+1]$ and $\eta_t$ is an independent Bernoulli random
variable with parameter $\frac{1+2\sigma}{2+2\sigma}$. We then have
$g'_{\sigma}(U)=\eta_t g'_{\sigma}(U_t)+(1-\eta_t)g'_{\sigma
}(V_t)$. Thus the law of $g'_{\sigma}(U)$ is the mixture of the law of
$g'_{\sigma}(U_t)$ and of the one of $g'_{\sigma}(V_t)$, with
respective weights $\frac{1+2\sigma}{2+2\sigma}$ and $1-\frac
{1+2\sigma}{2+2\sigma}$. Consequently
\[
\forall t\in[0,1],\ \forall x\in\R \qquad dP_{g'_{\sigma
}(U_t)}(x) \leq\frac{2+2\sigma}{1+2\sigma} \,dP_{g'_{\sigma}(U)}(x).
\]
The law of $Y_1(t)$ conditionally on $\#\{\tau_i\in[-\sigma,1+\sigma
]\}\ge1$ can be written as
\begin{eqnarray*}
&&dP_{Y_1(t)}(x)\\
&&\qquad=\frac{1}{1-e^{-\lambda(1+2\sigma)}} \sum
_{k=1}^{+\infty} e^{-\lambda(1+2\sigma)} \frac{(\lambda(1+2\sigma
))^k}{k!} \bigl(dP_{g'_{\sigma}(U_t)}\ast\cdots\ast dP_{g'_{\sigma}(U_t)}\bigr)(x).
\end{eqnarray*}
Thus, if we write $f_0= dP_{g'_{\sigma}(U)}$, we get
\begin{eqnarray*}
&&dP_{Y_1(t)}(x)\\
&&\qquad\leq\frac{1}{1-e^{-\lambda(1+2\sigma)}}
\sum_{k=1}^{+\infty} e^{-\lambda(1+2\sigma)} \frac{(\lambda(1+2\sigma
))^k}{k!} \biggl(\frac{2+2\sigma}{1+2\sigma}\biggr)^k (f_0 \ast
\cdots\ast f_0)(x) \\
&&\qquad= e^{\lambda}\frac{1-e^{-\lambda(2+2\sigma
)}}{1-e^{-\lambda(1+2\sigma)}} \tilde{f}_0(x),
\end{eqnarray*}
where $\tilde{f}_0(x)\,dx$ is a probability measure on $\R$. This shows
that we can take $c_1=e^{\lambda}\frac{1-e^{-\lambda(2+2\sigma
)}}{1-e^{-\lambda(1+2\sigma)}}$.

For $Y_2(t)$, we first notice that $Y_2(t)$ can be decomposed as the
sum of two independent random variables in the following way:
\begin{eqnarray*}
Y_2(t)&=& \sum_{\tau_i\in(-\infty,-1-\sigma+t]\cup[1+\sigma
+t,+\infty)} g'_{\sigma}(t-\tau_i) \\
&&{}+ \sum_{\tau_i\in(-\sigma
-1+t,-\sigma) \cup(1+\sigma,1+\sigma+t)} g'_{\sigma}(t-\tau_i) .
\end{eqnarray*}
The first random variable in the sum above has a law that does not
depend on $t$. For the second random variable, using the same trick as
above [{i.e.,} decompose here a uniform random variable on the interval
$(-1-\sigma,-\sigma)\cup(\sigma,1+\sigma)$ as a mixture with
weights $1/2$ and $1/2$ of two uniform random variables: one on
$(-1-\sigma,-1-\sigma+t)\cup(t+\sigma,1+\sigma)$, and the other
one on the rest], we obtain that $c_2=e^{\lambda}$.

And finally the bound on the expectation of the number of local extrema is
\begin{eqnarray*}
\rho(\lambda,\sigma) & \leq& \biggl( c_1 \frac{3\lambda(1+2\sigma
)}{1-e^{-\lambda(1+2\sigma)}} + c_2 \biggr) \bigl(1-e^{-\lambda
(1+2\sigma)}\bigr) + e^{-\lambda(1+2\sigma)} \\
& \leq& e^{\lambda}\frac{2+2\sigma}{1+2\sigma} \bigl(3\lambda
(1+2\sigma)\bigr) +e^{\lambda} =\bigl (3\lambda(2+2\sigma) +1\bigr) e^{\lambda} .
\end{eqnarray*}
\upqed\end{pf}

\subsection{Scaling properties}
An interesting property of the shot noise process with Gaussian kernel
is that we have two scale parameters: the intensity $\lambda$ of the
Poisson point process and the width $\sigma$ of the Gaussian kernel.
These two parameters are linked in the sense that changing one of them amounts
to change the other one in an appropriate way. These scaling properties
are described more precisely in the following lemma.

\begin{lemma}\label{scalinglem}
We have the following scaling properties for the process $X_{\lambda
,\sigma}$:
\begin{longlist}[1.]
\item[1.] Changing $\sigma$ and $\lambda$ in a proportional way: for
all $c>0$,
\[
\{ X_{{\lambda}/{c},c\sigma}(t) ; t\in\R\} \stackrel{\mathrm{f.d.d.}}{=}\biggl\{ \frac{1}{c} X_{\lambda,\sigma}\biggl(\frac{t}{c}\biggr); t\in\R \biggr\} .
\]
\item[2.] Increasing the width of the Gaussian kernel: for all $\sigma
_1$ and $\sigma_2$, we have
\[
\{ X_{\lambda,\sqrt{\sigma_1^2+\sigma_2^2}}(t) ; t\in\R\}
\stackrel{a.s.}{=} \{ (X_{\lambda,\sigma_1}\ast g_{\sigma_2})(t) ;
t\in\R\} .
\]
\item[3.] Increasing the intensity of the Poisson process: for all
$c>0$, we have
\[
\{ X_{\lambda\sqrt{1+c^2},\sigma}(t) ; t\in\R\} \stackrel{\mathrm{f.d.d.}}{=}
\bigl\{ \sqrt{1+c^2}\cdot(X_{\lambda,\sigma}\ast g_{c\sigma})\bigl(t\sqrt
{1+c^2}\bigr); t\in\R \bigr\} .
\]
\item[4.] The mean number $\rho(\lambda,\sigma)$ of local extrema of
$X_{\lambda,\sigma}$ per unit length satisfies
\[
\forall c > 0 \qquad  c\rho(\lambda, c\sigma)=\rho
(c\lambda,\sigma) .
\]
\end{longlist}
\end{lemma}

\begin{pf}
For the first property, let $\{\tau_i \}$ be a Poisson point process
of intensity $\lambda/c$ on the line. Then
\[
X_{{\lambda}/{c},c\sigma}(t)= \sum_i \frac{1}{c\sigma\sqrt
{2\pi}} e^{-(t-\tau_i)^2/2c^2\sigma^2} = \frac{1}{c} \sum_i
g_{\sigma}\biggl( \frac{t}{c} - \frac{\tau_i}{c} \biggr).
\]
Since the points $\{ \tau_i/c\} $ are now the points of a Poisson
process on intensity $\lambda$ on the line, we obtain the first
scaling property.
The second property comes simply from the fact that if $g_{\sigma_1}$
and $g_{\sigma_2}$ are two Gaussian kernels of respective width
$\sigma_1$ and $\sigma_2$, then their convolution is the Gaussian
kernel of width $\sqrt{\sigma_1^2+\sigma_2^2}$.
The third property is just a consequence of combining the first and
second properties.

For the fourth property, we first compute
\begin{eqnarray*}
X'_{\lambda,c\sigma}(t)&=& \frac{1}{c\sigma\sqrt{2\pi}}\sum_{\tau
_i} \frac{-(t-\tau_i)}{c^2\sigma^2} e^{-(t-\tau_i)^2/2c^2 \sigma
^2}\\
 &=& \frac{1}{c^2\sigma\sqrt{2\pi}}\sum_{\tau_i} \frac
{-(t/c-\tau_i/c)}{\sigma^2} e^{-(t/c-\tau_i/c)^2/2\sigma^2},
\end{eqnarray*}
where the $\{ \tau_i \}$ are the points of a Poisson point process of
intensity $\lambda$ on~$\R$. Then, since the $\{ \tau_i /c \}$ are
now the points of a Poisson point process of intensity $c \lambda$ on~$\R$,
we have that the expected number of points $t\in[0,c]$ such
that $X'_{\lambda,c\sigma}(t)=0$ [which, by definition, equals $c
\rho(\lambda,c \sigma)$], also equals the expected number of points
$t\in[0,1]$ such that $X^{'}_{c\lambda,\sigma}(t)=0$ [which is $\rho
(c \lambda,\sigma)$].
\end{pf}

\begin{figure}

\includegraphics{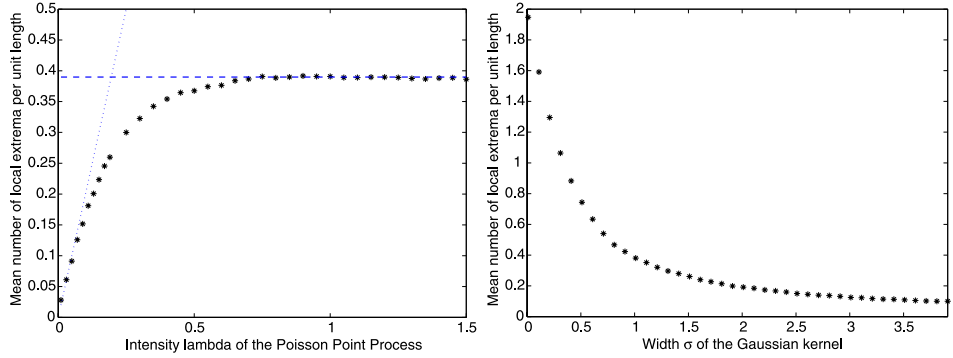}

\caption{On the left: empirical mean number of local extrema of
$X_{\lambda,\sigma}$ per unit length as a function of $\lambda$
(here $\sigma=1$ and we have taken the mean value from $50$ samples
on the interval $[-100,100])$. The horizontal dashed line is the
constant $\frac{1}{\pi}\sqrt{\frac{3}{2}}$ and the dotted line is
the map $\lambda\mapsto2\lambda$. On the right: empirical mean
number of local extrema of $X_{\lambda,\sigma}$ per unit length as a
function of $\sigma$ (here $\lambda=1$ and we have taken the mean
value from $10$ samples on the interval $[-100,100])$.}
\label{NbmaxlocLambdafig}
\end{figure}

To study how $\rho(\lambda,\sigma)$ varies when $\lambda$ and
$\sigma$ vary, we first can use the result on high intensity and
convergence to the crossings of a Gaussian process obtained in Theorem
\ref{convunifgausstheo}. Indeed, if the second moment of
$N_{X'_{\lambda,\sigma}}(0)$ is bounded by a polynomial in $\lambda
$, then
we will get
\[
\rho(\lambda,\sigma) \mathop{\longrightarrow}_{\lambda\to+\infty
} \frac{1}{\sigma\pi}\sqrt{\frac{3}{2}} .
\]
And thanks to the scaling properties, this also will imply that $\rho
(\lambda,\sigma)$ is equivalent to $\frac{1}{\sigma\pi}\sqrt
{\frac{3}{2}}$ as $\sigma$ goes to $+\infty$.
These two facts have been empirically checked and are illustrated on
Figure~\ref{NbmaxlocLambdafig}. Now, notice that we can also observe
on the left-hand figure another regime when $\lambda$ is small.
Indeed, $\rho(\lambda,\sigma)$ seems to be almost linear for small
values of $\lambda$. Notice also on the right-hand figure that $\rho
(\lambda,\sigma)$ seems to be\vadjust{\goodbreak} a decreasing function of $\sigma$
(this then indicates that, as time goes by, the temperature on the line
fluctuates less and less). The study of these two facts is the aim of
the next section.

\subsection{Heat equation and local extrema}

In this subsection we assume first that $\lambda>0$ is fixed.
As we already mentioned it in the introduction of Section~\ref{Gaussiankernel},
one of the main features of the shot noise process $X_{\lambda,\sigma
}$ is
that it can be seen in a dynamic way, which means that we can study how it
evolves as the width $\sigma$ of the Gaussian kernel changes and
consider it as a random field indexed by the variable $(\sigma,t)$.
Then, the main tool is the heat equation which is satisfied by the
Gaussian kernel
%
\begin{eqnarray}\label{heateq}
&&\forall\sigma>0,\   \forall t\in\R
\nonumber
\\[-8pt]
\\[-8pt]
\nonumber
&&\qquad \frac{\partial
g_{\sigma}}{\partial\sigma} (t) = {\sigma}   g''_{\sigma}(t)
\quad\mbox{and also consequently}\quad \frac{\partial g'_{\sigma}}{\partial
\sigma} (t) = {\sigma}   g^{(3)}_{\sigma}(t) .
\end{eqnarray}
Since the Gaussian kernel $g_\sigma$ is a very smooth function, both
in $\sigma>0$ and $t\in\R$, by the same type of proof as the ones in
Proposition~\ref{diff}, we have that $(\sigma,t)\mapsto X_{\lambda
,\sigma}(t)$ is almost surely and mean square smooth on $(0,+\infty
)\times\R$ with
%
\begin{eqnarray}\label{heatforXeq}
\frac{\partial X_{\lambda,\sigma}}{\partial\sigma} (t)& =& \sum_i
\frac{\partial g_{\sigma}}{\partial\sigma} (t-\tau_i) = {\sigma}
  X''_{\lambda,\sigma}(t)\quad \mbox{and also}
  \nonumber
  \\[-8pt]
  \\[-8pt]
  \nonumber
  \frac{\partial
X'_{\lambda,\sigma}}{\partial\sigma} (t) &=& {\sigma}
X^{(3)}_{\lambda,\sigma}(t) .
\end{eqnarray}
We will see in the following that this equation will be of great
interest to study the crossings of $X_{\lambda,\sigma}$.\vadjust{\goodbreak}

The convolution of a real function defined on $\R$ with a Gaussian
kernel of increasing width $\sigma$ (which amounts to apply the heat
equation) is a very common smoothing technique in signal processing.
One of its main properties is generally formulated by the wide-spread
idea that ``Gaussian convolution on $\R$ cannot create new extrema''
(and it is in some sense the only kernel that has this property; see
\cite{YuillePoggio}). This has been studied (together with its
extension in higher dimension) for applications in image processing by
Lindeberg~\cite{Lindeberg}, and also by other authors (e.g.,
to study mixtures of Gaussian distributions as in \cite
{CarreiraWilliams1} and~\cite{CarreiraWilliams2}).
However, in most cases, the correct mathematical framework for the
validity of this property is not exactly stated. Thus we start here
with a lemma giving the conditions under which one can obtain
properties for the zero-crossings of a function solution of the heat equation.
The result, which proof is postponed to the \hyperref[Appendix]{Appendix}, is stated under a
general form for a function $h$ in
the two variables $\sigma$ and $t$. But we have to keep in mind that
we will want to apply this to $h(\sigma,t)=X'_{\lambda,\sigma}(t)$
to follow the local extrema of the shot noise process when $\sigma$ evolves.

\begin{lemma}\label{PathOfZeroslem}
Let $\sigma_0>0$ and $(\sigma,t)\mapsto h(\sigma,t)$ be a ${\mathcal
C}^2$ function defined on $(0,\sigma_0]\times[a,b]$, which satisfies
the heat equation
\[
\forall(\sigma,t)\in(0,\sigma_0]\times\R  \qquad \frac
{\partial h}{\partial\sigma}(\sigma,t) = {\sigma}   \frac
{\partial^2 h}{\partial t^2}(\sigma,t) .
\]
We assume that:
\begin{longlist}[(a)]
\item[(a)] there are no $t\in[a,b]$ such that $h(\sigma_0,t)=0$
and $\frac{\partial h}{\partial t}(\sigma_0,t)=0$,
\item[(b)] there are no $(\sigma,t)\in(0,\sigma_0]\times[a,b]$
such that $h(\sigma,t)=0$ and \mbox{$\nabla h(\sigma,t)=0$}.
\end{longlist}
Then we have the following properties for the zero-crossings of $h$:
\begin{longlist}[(iii)]
\item[(i)] Global curves: If $t_0\in(a,b)$ is such that $h(\sigma
_0,t_0)=0$, there exists $\sigma_0^-<\sigma_0$ and a maximal
continuous path $\sigma\mapsto\Gamma_{t_0}(\sigma)$ defined on
$(\sigma_0^-,\sigma_0]$ such that $\Gamma_{t_0}(\sigma_0)=t_0$ and
for all $\sigma\in(\sigma_0^-,\sigma_0]$ we have $h(\sigma,\Gamma
_{t_0}(\sigma))=0$. Moreover, if $\Gamma_{t_0}(\sigma)$ stays within
some compact set of $\R$ for all $\sigma$, then $\sigma_0^-=0$.
\item[(ii)] Nonintersecting curves: If $\widetilde{t_0}\neq t_0$
is another point in $(a,b)$ such that $h(\sigma_0,\widetilde
{t_0})=0$, then for all $\sigma\in(0,\sigma_0]$ we have $\Gamma
_{t_0}(\sigma)\neq\Gamma_{\widetilde{t_0}}(\sigma)$.
\item[(iii)] Local description of the curves: If $(\sigma
_1,t_1)\in(0,\sigma_0]\times\R$ is such that $h(\sigma_1,t_1)=0$
then there exist a ${\mathcal C}^1$ function $\eta$ defined on a
neighborhood of $\sigma_1$ and such that $h(\sigma,\eta(\sigma))=0$
in this neighborhood of $\sigma_1$, or a ${\mathcal C}^1$ function
$\xi$ defined on a neighborhood of $t_1$ and such that $h(\xi
(t),t)=0$ in this neighborhood of~$t_1$, and moreover, if $\xi
'(t_1)=0$, then $\xi''(t_1)<0$ (it is a local maximum).
\end{longlist}
\end{lemma}

\begin{figure}

\includegraphics{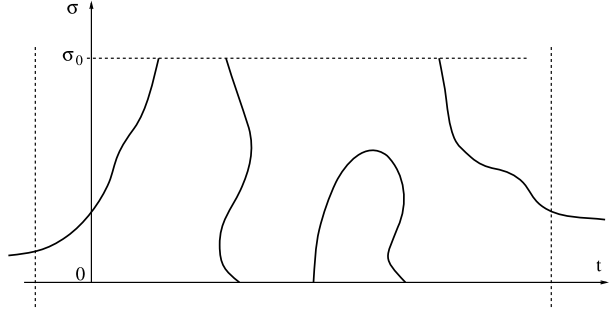}

\caption{Curves of $h(\sigma,t)=0$ for some $h$ satisfying the heat
equation, in the $(t,\sigma)$ domain; here $t$ is along the horizontal
axis and $\sigma$ is along the vertical one. According to Lemma \protect\ref
{PathOfZeroslem}, the zeros-crossings of $h$ are a set of
nonintersecting curves, that are locally else functions of $\sigma$ or
functions of $t$ with no local minima.}
\label{PathOfZerosfig}
\end{figure}

The properties stated in Lemma~\ref{PathOfZeroslem} are illustrated on
Figure~\ref{PathOfZerosfig}, where the different types of curves
formed by the set of points $\{(t,\sigma)\in\R^2; h(\sigma,t)=0\}$
are shown for some $h$ satisfying the heat equation.

Let us consider again the shot noise process $X_{\lambda,\sigma}$. We
now give the main result for the number of local extrema of $X_{\lambda
,\sigma}$ as a function of $\sigma$. The intensity $\lambda$ is
assumed to be fixed.

\begin{theorem}\label{SigmaVarth}
Let $\sigma_0>0$ and $a\le b$. Then,
\[
\Prob\bigl( \exists(\sigma,t)\in(0,\sigma_0)\times[a,b] \mbox{ such
that } X'_{\lambda,\sigma}(t)=0 \mbox{
and } \nabla X'_{\lambda,\sigma}(t)=0 \bigr) = 0 .
\]
Moreover, if we assume that for all $0<\sigma_1<\sigma_0$
\[
\E\bigl( \#\{\sigma\in[\sigma_1,\sigma_0] \mbox{ such that }
X'_{\lambda,\sigma}(0)=0 \}\bigr) < +\infty,
\]
then the function $\sigma\mapsto\rho(\lambda,\sigma)$, which gives
the mean number of local extrema of $X_{\lambda,\sigma}$ per unit
length, is
decreasing and it has the limit $2\lambda$ as $\sigma$ goes to $0$.
\end{theorem}

\begin{pf}
Let us denote $Y(\sigma,t):=X'_{\lambda,\sigma}(t)$ for all $(\sigma
,t)\in(0,+\infty)\times\R$. We first check that the assumptions (a)
and (b) of Lemma~\ref{PathOfZeroslem} are satisfied almost surely for
$Y$. Assumption (a) is already given by Proposition~\ref{bound19}. For
assumption (b), we first notice that since $Y(\sigma,t)$ satisfies the
heat equation, we have
\begin{eqnarray*}
&&\{ Y({\sigma},t)=0 \mbox{ and } \nabla Y(\sigma,t)=0 \}
\\
&&\qquad=\{ Y({\sigma},t)=0 \mbox{
and } Y'(\sigma,t)=0 \mbox{ and } Y''({\sigma},t)=0 \}.
\end{eqnarray*}
Then a slight modification of the proof of Proposition \ref
{NoTangencyprop}, using the second-order Taylor formula in \eqref{Taylor1},
allows us to conclude that
$\Prob( \exists(\sigma,t)\in(0,\sigma_0)\times[a,b] \mbox{ such
that } Y({\sigma},t)=0 \mbox{ and } \nabla Y(\sigma,t)=0 ) = 0$,
using the same integrability bound for the characteristic function of
$Y(\sigma,t)$ as the one obtained in the proof of Proposition \ref
{bound19} [and considering first $(\sigma,t)\in(\sigma_1,\sigma_0)
$ for $\sigma_1>0$, and conditioning by $\{\gamma_T\geq3\}$].
This also proves the first part of the theorem.

Let $0<\sigma_1<\sigma_0$ be fixed. By assumption, we have $ \E( \#\{
\sigma\in[\sigma_1,\sigma_0]$ such that $X'_{\lambda
,\sigma}(0)=0 \}) < +\infty.$
Notice that by stationarity this expected value is independent of the
value of $t$ (taken as $0$ above).
Let $T>0$ and let us consider the zeros of $Y(\sigma,t)= X'_{\lambda
,\sigma}(t)$ for $(\sigma,t)\in[\sigma_1,\sigma_0]\times[0,T]$.
Let $t_0\in[0,T]$ be such that $Y(\sigma_0,t_0)=0$. By Lemma \ref
{PathOfZeroslem}, there is a continuous path $\sigma\mapsto\Gamma
_{t_0}(\sigma)$ that will ``cross the left or right boundary of the
domain,'' that is, be such that there exists $\sigma\in$ such that
$\Gamma_{t_0}(\sigma)=0 \mbox{ or } T$, or will be defined until
$\sigma_1$ and such that $\Gamma_{t_0}(\sigma_1)\in[0,T]$. We thus
have 
\[
\rho(\sigma_0,[0,T]) \leq2 \E\bigl( \#\{\sigma\in[\sigma_1,\sigma_0]
\mbox{ such that } X'_{\lambda,\sigma}(0)=0\} \bigr) + \rho(\sigma
_1,[0,T]).
\]
Dividing both sides by $T$ and letting $T$ go to infinity then shows
that $\rho(\sigma_0) \leq\rho(\sigma_1)$. Thus the function
$\sigma\mapsto\rho(\lambda,\sigma)$ is decreasing.

To find the limit of $\rho(\lambda,\sigma)$ as $\sigma$ goes to $0$
(that exists thanks to the bound of Proposition~\ref{bound19}),
instead of looking at the local extrema of $X_{\lambda,\sigma}$ in
$[0,1]$, we will only look at the local maxima (which are the
down-crossings of $0$ by the derivative) in $[0,1]$. Let
$D_{X'_{\lambda,\sigma}}(0,[0,1])$ be the random variable that counts
these local maxima, and let $\rho^-(\lambda,\sigma)=\E
(D_{X'_{\lambda,\sigma}}(0,[0,1]))$. By stationarity of $X_{\lambda
,\sigma}(t)$ and because between any two local maxima, there is a
local minima, we have that $\rho^-(\lambda,\sigma)= \frac{1}{2}\rho
(\lambda,\sigma)$.
Now, we introduce ``barriers'' in the following way: let $E_{\sigma
_0}$ be the event ``there are no points of the Poisson point process in
the intervals $[-2\sigma_0,2\sigma_0]$ and $[1-2\sigma_0,1+2\sigma_0]$.''
If we assume that $E_{\sigma_0}$ holds, then $X''_{\lambda,\sigma
}(t)>0$ for all $t$ in $[-\sigma_0,\sigma_0]\cup[1-\sigma
_0,1+\sigma_0]$ and all $\sigma\leq\sigma_0$, and therefore there
are no local maxima of $X_{\lambda,\sigma}$ in these intervals.
Then by Lemma~\ref{PathOfZeroslem}, we can follow all the local
maxima of $X_{\lambda,\sigma}$ in $[0,1]$ from $\sigma=\sigma_0$
down to $\sigma=0$. Thus $\sigma\mapsto D_{X'_{\lambda,\sigma
}}(0,[0,1])\ind_{E_{\sigma_0}}$ is a decreasing function of $\sigma$
for $\sigma\leq\sigma_0$.
Moreover, we can also check that the set of local maxima of $X_{\lambda
,\sigma}(t)$ in $[0,1]$ converges, as $\sigma$ goes to $0$, to the
set of points of the Poisson process in $[0,1]$. This implies, in
particular, that $D_{X'_{\lambda,\sigma}}(0,[0,1])$ goes to $\#\{\tau
_i\in[0,1] \}$ as $\sigma$ goes to $0$. Thus by monotone convergence,
it implies that $\rho^-(\lambda,\sigma| E_{\sigma_0})$ goes to $\E
(\#\{\tau_i\in[0,1] \}| E_{\sigma_0})$. Since the sequence of events
$E_{\sigma_0}$ is an increasing sequence of events as $\sigma_0$
decreases to $0$, we finally get
\[
\lim_{\sigma\to0}\rho^-(\lambda,\sigma) =
\lim_{\sigma_0 \to0} \E(\#\{\tau_i\in[0,1] \}| E_{\sigma_0}) =
\E(\#\{\tau_i\in[0,1] \}) = \lambda.
\]
\upqed\end{pf}

\begin{figure}

\includegraphics{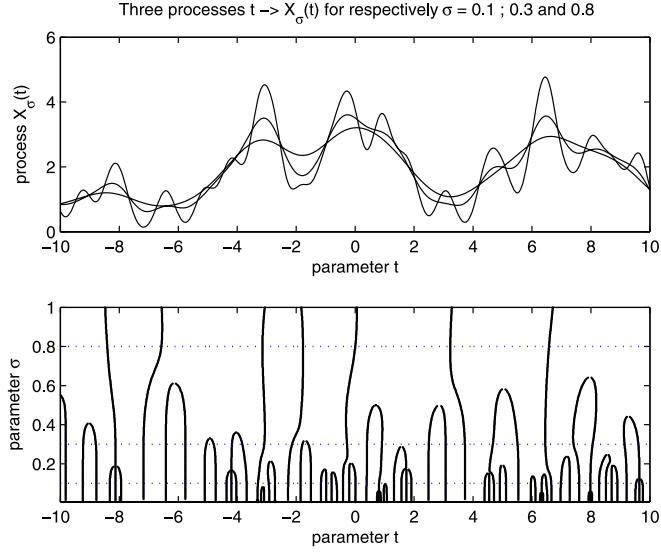}

\caption{Top: three processes $t\mapsto X_{\lambda,\sigma}(t)$
obtained from the same Poisson point process of intensity $\lambda=2$
and for a Gaussian kernel of respective width $\sigma=0.1; 0.3$ and
$0.8$. Bottom: evolution of the local extrema of $t\mapsto X_{\lambda
,\sigma}(t)$ as $\sigma$ goes from $0$ to $1$. The three values
$\sigma=0.1; 0.3$ and $0.8$ are plotted as dotted line. They indicate
the local extrema of the three processes above.}
\label{FollowLocExtrfig}
\end{figure}

Thus, under the assumption that
$ \E( \#\{\sigma\in[\sigma_1,\sigma_0] \mbox{ such that }
X'_{\lambda,\sigma}(0)=0 \}) < +\infty$ for all $0<\sigma_1<\sigma
_0$, Theorem~\ref{SigmaVarth} asserts that the function $\sigma
\mapsto\rho(\lambda,\sigma)$ is a decreasing function with limit
$2\lambda$ when $\sigma\to0$. This fact was empirically observed on
Figure~\ref{NbmaxlocLambdafig}, and is also illustrated on Figure
\ref{FollowLocExtrfig} where we ``follow'' the local extrema as
$\sigma$ evolves. Now, these properties can be translated, using the
scaling relations of Lemma~\ref{scalinglem}, into the following
properties on $\lambda\mapsto\rho(\lambda,\sigma)$:
\[
\forall c\geq1 ,   \rho(c\lambda,\sigma)\leq c \rho
(\lambda,\sigma) ;  \qquad \rho(\lambda,\sigma)\leq
2\lambda  \quad \mbox{and}\quad  \frac{\rho
(\lambda,\sigma)}{2\lambda}\mathop{\longrightarrow}_{\lambda\to0}1 .
\]
This shows the second asymptotic linear regime observed for small
values of the intensity $\lambda$.

\begin{appendix}
\section*{Appendix}\label{Appendix}
\subsection{Stationary phase estimate for oscillatory integrals}
%
\begin{proposition}[(Stationary phase estimate for oscillatory integrals)]\label{PhaseStatiolem}
Let $a<b$ and let $\varphi$ be a function of class $\mathcal{C}^2$
defined on $[a,b]$. Assume that $\varphi'$ and $\varphi''$ cannot
simultaneously vanish on $[a,b]$ and denote $  m=\break \min
_{s\in[a,b]}\sqrt{\varphi'(s)^2+\varphi''(s)^2} >0$. Let us also
assume that $n_0=\#\{s\in[a,b] \mbox{ s.t.}\break \varphi''(s)=0 \}
<+\infty$. Then
\[
\forall u\in\R\mbox{ s.t. } |u|>\frac{1}{m} \qquad
\biggl| \int_a^b e^{iu\varphi(s)}  \,ds \biggr| \leq\frac{8 \sqrt{2}
(2n_0+1)}{\sqrt{m|u|}} .
\]

Now, let $\varphi_1$ and $\varphi_2$ be two functions of class
$\mathcal{C}^3$ defined on $[a,b]$. Assume that the derivatives of
these functions are linearly independent, in the sense that for all
$s\in[a,b]$, the matrix $\Phi(s)=\bigl(
{\varphi_1'(s) \atop \varphi_1''(s)}\enskip{\varphi_2'(s) \atop
\varphi_2''(s)}\bigr)$ is invertible. Denote $  m=\min_{s\in
[a,b]}\Vert \Phi(s)^{-1}\Vert ^{-1}>0$, where $\|\cdot\|$ is
the matricial norm induced by the Euclidean one. Assume, moreover, that
there exists $n_0<+\infty$ such that $\#\{s\in[a,b] \mbox{ s.t. }
{\operatorname{det}}(\Phi'(s)) =0 \} \leq n_0$, where $\Phi'(s)=\bigl(
{ \varphi_1''(s) \atop  \varphi
_1^{(3)}(s)}\enskip{\varphi_2''(s) \atop \varphi_2^{(3)}(s)}
\bigr)$. Then
\[
\forall(u,v)\in\R^2 \mbox{ s.t. } \sqrt{u^2+v^2}>\frac{1}{m} \qquad
   \biggl| \int_a^b e^{iu\varphi_1(s)+iv\varphi_2(s)}
\,ds \biggr| \leq\frac{8 \sqrt{2} (2n_0+3)}{\sqrt{m\sqrt{u^2+v^2}}} .
\]
\end{proposition}

\begin{pf}
For the first part of the proposition, by assumption, $[a,b]$ is the
union of the three compact sets
\begin{eqnarray*}
&&\{s\in[a,b]; |\varphi''|\ge m/2\},\qquad \{s\in[a,b];
|\varphi'|\ge m/2 \mbox{ and } \varphi'' \ge0 \} \quad\mbox{and }\\
&&\qquad\{s\in[a,b]; |\varphi'|\ge m/2 \mbox{ and } \varphi'' \le0
\}.
\end{eqnarray*}
Therefore there exists $1\leq n \leq2n_0+1$ and a subdivision
$(a_i)_{0\le i\le n}$ of $[a,b]$ such that $[a_{i-1},a_i]$ is included
in one of the previous subsets for any $1\le i\le n$.
If $[a_{i-1}, a_i]\subset\{s\in[a,b]; |\varphi''(s)|\ge
m/2\}$, according to~\cite{Stein}, Proposition 2, page~332,
\[
\biggl|\int_{a_{i-1}}^{a_i} e^{iu \varphi(s)}  \,ds\biggr| =
\biggl|\int_{a_{i-1}}^{a_i} e^{iu(m/2) (2 \varphi(s)/m)}  \,ds\biggr|\le
8\frac{\sqrt{2}}{\sqrt{m|u|}},
\]
otherwise,
\[
\biggl|\int_{a_{i-1}}^{a_i} e^{iu \varphi(s)}  \,ds\biggr| \le\frac
{6}{m |u|}.
\]
The result follows from summing up these $n$ integrals.

For the second part of the proposition, we use polar coordinates, and
write $(u,v)=(r\cos\theta,r\sin\theta)$. For $\theta\in[0,2\pi
)$, let $\varphi_{\theta}$ be the function defined on $[a,b]$ by
$\varphi_{\theta}(s)=\varphi_1(s) \cos\theta+\varphi_2(s) \sin
\theta$. Then $
{\varphi'_{\theta}(s) \choose \varphi''_{\theta}(s)} = \Phi(s)
{\cos\theta
\choose  \sin\theta}$, and thus $1=\Vert\Phi(s)^{-1}
{\varphi'_{\theta}(s) \choose \varphi''_{\theta}(s)} \Vert$. This implies that for all $s\in[a,b]$, $\sqrt
{\varphi'_{\theta}(s)^2+\varphi''_{\theta}(s)^2} \geq1/\Vert
\Phi(s)^{-1}\Vert \geq m$. Moreover, thanks to Rolle's theorem,
the number of points $s\in[a,b]$ such that $\varphi''_{\theta}(s)=0$
is bounded by one plus the number of $s\in[a,b]$ such that $\varphi
''_{1}(s)\varphi'''_{2}(s) - \varphi'''_{1}(s)\varphi''_{2}(s) =0$,
that is, by $1+n_0$.
Thus, we can apply the result of the first part of the proposition to
each function $\varphi_{\theta}$ and the obtained bound will depend
only on $m$, $n_0$ and $r=\sqrt{u^2+v^2}$.
\end{pf}

\subsection{\texorpdfstring{Proof of Proposition \protect\ref{BoundsConvprop}}
{Proof of Proposition 8}}

For $k\geq0$ and $l\geq0$ integers, let us denote $m_{kl}=\int
|g(s)|^k |g'(s)|^l   \,ds$. We will also simply denote $m_0=m_{20}=\int
g(s)^2  \,ds$ and $m_2=m_{02}=\int g'(s)^2 \,ds$.

Let $\psi_{{\lambda}}(u,v)$ denote the joint characteristic function
of $(Z_{\lambda}(t),Z'_{\lambda}(t))$, then
\begin{eqnarray*}
\psi_{{\lambda}}(u,v) &=& \E\bigl( e^{i({u}/{\sqrt{\lambda}})X_{\lambda}+i({v}/{\sqrt{\lambda}})X'_{\lambda}}\bigr) e^{-iu\sqrt
{\lambda}\int g}\\
& =& \exp\biggl(\lambda\int_{\R} \biggl(e^{i
({u}/{\sqrt{\lambda}})g(s) +i({v}/{\sqrt{\lambda}})g'(s)} - 1
-i\frac{u}{\sqrt{\lambda}} g(s) \biggr)   \,ds\biggr).
\end{eqnarray*}
We now use the fact $\int g' =0$, and we thus have $\psi_{{\lambda
}}(u,v) = \exp(H_{\lambda}(u,v))$ where
\[
H_{\lambda}(u,v)= \lambda\int_{\R} \biggl(e^{i({u}/{\sqrt
{\lambda}})g(s) +i({v}/{\sqrt{\lambda}})g'(s)} - 1 -i
\frac{u}{\sqrt{\lambda}} g(s)-i\frac{v}{\sqrt{\lambda}} g'(s) \biggr)
  \,ds .
\]
We need to notice that
\[
\forall(u,v)\in\R^2  \qquad |\psi_{{\lambda}}(u,v)| =
|\exp(H_{\lambda}(u,v))| = |\E( e^{iuZ_{\lambda}+ivZ'_{\lambda}})|
\leq1 .
\]
In the following, we will also need these simple bounds:
%
\begin{eqnarray}\label{simplebounds}
\forall x\in\R \qquad  \biggl|e^{ix} -1 -ix + \frac{x^2}{2} \biggr|
&\leq&\frac{|x|^3}{3!}  \quad \mbox{and}
\nonumber
\\[-8pt]
\\[-8pt]
\nonumber
\forall z\in\C\qquad   |e^z - 1|&\leq&|z| e^{|z|} .
\end{eqnarray}

We first estimate $H_{\lambda}(u,0)$. We have
\[
H_{\lambda}(u,0)=\lambda\int\biggl(e^{i({u}/{\sqrt{\lambda}})g(s)} -
1 -i\frac{u}{\sqrt{\lambda}} g(s)\biggr)   \,ds = -\frac{1}{2} u^2 m_0 +
K_{\lambda}(u),
\]
where $K_{\lambda}(u)=\lambda\int(e^{i({u}/{\sqrt{\lambda
}})g(s)} - 1 -i\frac{u}{\sqrt{\lambda}} g(s) + \frac{1}{2} \frac
{u^2}{\lambda} g^2(s))   \,ds$. Then, thanks to the simple bounds (\ref
{simplebounds}), we get
\begin{eqnarray*}
|K_{\lambda}(u)|&\leq&\frac{|u|^3 m_{30}}{6\sqrt{\lambda}}\qquad \mbox{and consequently} \\
\bigl|e^{H_{\lambda}(u,0)} - e^{-({1}/{2}) u^2 m_0}\bigr|
&\leq&\frac{|u|^3 m_{30}}{6\sqrt{\lambda}} e^{-({1}/{2}) u^2 m_0 }
e^{{|u|^3 m_{30}}/{(6\sqrt{\lambda})}} .
\end{eqnarray*}
We then estimate $H_{\lambda}(u,v) - H_{\lambda}(u,0)$,
\begin{eqnarray*}
H_{\lambda}(u,v) - H_{\lambda}(u,0) & = & \lambda\int\bigl(e^{i
({u}/{\sqrt{\lambda}})g(s) +i({v}/{\sqrt{\lambda}})g'(s)} -
e^{i({u}/{\sqrt{\lambda}})g(s)}\bigr)   \,ds \\
& = & \lambda\int e^{i({u}/{\sqrt{\lambda}})g(s)}\biggl(e^{i
({v}/{\sqrt{\lambda}})g'(s)} -1 -i\frac{v}{\sqrt{\lambda}} g'(s)\biggr)
\,ds \\
& = & -\frac{v^2}{2} \int g'(s)^2 e^{i({u}/{\sqrt{\lambda}})g(s)}
  \,ds + F_{\lambda}(u,v) ,
\end{eqnarray*}
where $F_{\lambda}(u,v)=\lambda\int e^{i({u}/{\sqrt{\lambda
}})g(s)}(e^{i({v}/{\sqrt{\lambda}})g'(s)} -1 -i\frac{v}{\sqrt
{\lambda}} g'(s)+ \frac{v^2}{2\lambda} g'(s)^2)   \,ds $. And again,
thanks to the simple bounds (\ref{simplebounds}),\vadjust{\goodbreak} we get $|F_{\lambda
}(u,v) |\leq\frac{|v|^3 m_{03}}{6\sqrt{\lambda}}$.
This implies that
\begin{eqnarray*}
&&\bigl|e^{H_{\lambda}(u,v) - H_{\lambda}(u,0)} - e^{-{v^2}/{2}
\int g'(s)^2
e^{i({u}/{\sqrt{\lambda}})g(s)}   \,ds}\bigr| \\
&&\qquad \leq
\bigl|e^{-{v^2}/{2} \int
g'(s)^2 e^{i({u}/{\sqrt{\lambda}})g(s)}   \,ds} \bigr| \cdot
\bigl|e^{F_{\lambda}(u,v)} -1 \bigr| \\
&&\qquad \leq  \frac{|v|^3 m_{03}}{6\sqrt{\lambda}} e^{ -{v^2}/{2}
\int g'(s)^2
\cos(({u}/{\sqrt{\lambda}})g(s))   \,ds +{|v|^3
m_{03}}/({6\sqrt{\lambda}})} .
\end{eqnarray*}

Let us now compute $\widehat{C_{Z_{\lambda}}}(u,[0,1])$. By
Proposition~\ref{Crossingspropgenerale}, we know that
\[
-\pi\widehat{C_{Z_{\lambda}}}(u,[0,1]) = \int_{0}^{+\infty} \frac
{1}{v^2}\bigl(\psi_{{\lambda}}(u,v) + \psi_{{\lambda}}(u,-v) -
2\psi_{{\lambda}}(u,0) \bigr)   \,dv .
\]
Let $V>0$ be a real number. We split the integral above in two parts,
and write it as the sum of the integral between $0$ and $V$, and of the
integral between $V$ and $+\infty$. Since for all $(u,v)$, we have
$|\psi_{{\lambda}}(u,v)|\leq1$, we get
\[
\biggl| \int_{V}^{+\infty} \frac{1}{v^2}\bigl(\psi_{{\lambda
}}(u,v) + \psi_{{\lambda}}(u,-v) - 2\psi_{{\lambda}}(u,0) \bigr)
  \,dv \biggr| \leq4 \int_{V}^{+\infty} \frac{1}{v^2}  \,dv = \frac
{4}{V} .
\]
On the other hand, let $I_V(u)$ denote the integral between $0$ and
$V$. We have
\[
I_V(u) = \int_{0}^{V} \frac{1}{v^2}e^{H_{\lambda}(u,0)}\bigl(
e^{H_{\lambda}(u,v)-H_{\lambda}(u,0)} +e^{H_{\lambda
}(u,-v)-H_{\lambda}(u,0)} -2 \bigr)   \,dv .
\]
We then decompose this into
\begin{eqnarray*}
I_V(u) & = & \int_{0}^{V} \frac{1}{v^2}e^{H_{\lambda}(u,0)}\bigl(
e^{H_{\lambda}(u,v)-H_{\lambda}(u,0)}
+e^{H_{\lambda}(u,-v)-H_{\lambda}(u,0)} \\
&&\hspace*{102pt}{}- 2 e^{-{v^2}/{2} \int g'(s)^2
e^{i({u}/{\sqrt{\lambda}})g(s)}   \,ds}\bigr)  \,dv \\
& & {}+ \int_{0}^{V} \frac{1}{v^2}e^{H_{\lambda}(u,0)}\bigl(2
e^{-{v^2}/{2}
\int g'(s)^2 e^{i({u}/{\sqrt{\lambda}})g(s)}   \,ds} - 2
e^{-({v^2}/{2})m_2}\bigr)  \,dv \\
& & {}+ \int_{0}^{V}
\frac{1}{v^2}\bigl(e^{H_{\lambda}(u,0)}- e^{-({1}/{2}) u^2
m_0}+e^{-({1}/{2})
u^2 m_0} \bigr)\bigl( 2 e^{-({v^2}/{2})m_2}-2 \bigr)  \,dv .
\end{eqnarray*}
Using the bounds we computed above, we get that
\begin{eqnarray*}
&&\biggl|I_V(u)- 2e^{-({1}/{2}) u^2 m_0} \int_0^V\frac{e^{-
({v^2}/{2})m_2}-1}{v^2}
  \,dv \biggr|\\
   & & \qquad\leq2\int_{0}^{V} \frac{v m_{03}}{6\sqrt{\lambda
}} e^{
-{v^2}/{2}
\int g'(s)^2 \cos(({u}/{\sqrt{\lambda}})g(s))   \,ds +{|v|^3 m_{03}}/({6\sqrt{\lambda}})}  \,dv \\
& &\qquad\quad{} + 2
\int_{0}^{V}\frac{1}{v^2}\bigl|e^{-{v^2}/{2} \int g'(s)^2
e^{i({u}/{\sqrt{\lambda}})g(s)}   \,ds} - 2 e^{-
({v^2}/{2})m_2}\bigr|  \,dv \\
& &\qquad\quad{} + 2 \bigl|e^{H_{\lambda}(u,0)}- e^{-({1}/{2}) u^2 m_0}\bigr|
\int_0^V\frac{1-e^{-({v^2}/{2})m_2}}{v^2}   \,dv .
\end{eqnarray*}
Let $J^{(n)}_V(u)$, for $n=1,2,3,$ respectively, denote the three terms above.
To give an upper bound for $J^{(1)}_V(u)$, we will need the following
basic inequality: $\forall x\in\R$, $\cos(x) \geq1-\frac{x^2}{2}$.
This gives us the bound
\[
J^{(1)}_V(u)\leq2\int_{0}^{V} \frac{v m_{03}}{6\sqrt{\lambda}} e^{
-{v^2m_2}/{2}+({v^2}/{2})({u^2m_{22}}/{(2\lambda)}) +
{|v|^3 m_{03}}/{(6\sqrt{\lambda}})}  \,dv .
\]
For the second term, we use
\begin{eqnarray*}
&&\bigl|e^{-{v^2}/{2} \int g'(s)^2 e^{i({u}/{\sqrt{\lambda}})g(s)}   \,ds} - 2
e^{-({v^2}/{2})m_2}\bigr| \\
&&\qquad \leq e^{-({v^2}/{2})m_2} \bigl|
e^{-({v^2}/{2}) \int
g'(s)^2 (e^{i({u}/{\sqrt{\lambda}})g(s)}-1)   \,ds} -1 \bigr| \\
&&\qquad \leq e^{-({v^2}/{2})m_2}\biggl|\frac{v^2}{2} \int g'(s)^2
\bigl(e^{i({u}/{\sqrt{\lambda}})g(s)}-1\bigr)   \,ds\biggr| e^{|
({v^2}/{2}) \int g'(s)^2
(e^{i({u}/{\sqrt{\lambda}})g(s)}-1)   \,ds|} .
\end{eqnarray*}
But $|\int g'(s)^2 (e^{i({u}/{\sqrt{\lambda}})g(s)}-1)
\,ds|\leq\int g'(s)^2\frac{|u|}{\sqrt{\lambda}}g(s)   \,ds =
\frac{|u|}{\sqrt{\lambda}}m_{12}$ and thus
\[
J^{(2)}_V(u)\leq\frac{|u|}{\sqrt{\lambda}}m_{12}\int
_{0}^{V}e^{-({v^2}/{2})m_2 + ({v^2}/{2})({|u|}/{\sqrt
{\lambda}})m_{12}}   \,dv.
\]
For the third term, we use an integration by parts to obtain that
\[
\int_0^V\frac{1-e^{-({v^2}/{2})m_2}}{v^2}   \,dv = \frac{e^{-
({V^2}/{2})m_2}-1}{V} + \int_0^V m_2 e^{-({v^2}/{2})m_2}   \,dv \leq
\frac{1}{2}\sqrt{2\pi m_2},
\]
which gives
\[
J^{(3)}_V(u)\leq\sqrt{2\pi m_2}\frac{|u|^3 m_{30}}{6\sqrt{\lambda
}} e^{-({1}/{2}) u^2 m_0 +{|u|^3 m_{30}}/{(6\sqrt{\lambda})}} .
\]
Moreover, we also have
\begin{eqnarray*}
\biggl|2 \int_0^V\frac{1-e^{-({v^2}/{2})m_2}}{v^2}   \,dv -\sqrt
{2\pi m_2}\biggr| &\leq&\frac{1-e^{-({V^2}/{2})m_2}}{V} + \int
_V^{+\infty} m_2 e^{-({v^2}/{2})m_2}   \,dv\\
 &\leq&\frac{2}{V}.
\end{eqnarray*}
The partial conclusion of all these estimates is that
\begin{eqnarray*}
&&\bigl|\pi\widehat{C_{Z_{\lambda}}}(u,[0,1]) - \sqrt{2\pi m_2}
e^{-m_0 u^2/2}\bigr| \\
&&\qquad\leq\frac{4}{V}+ \frac{2e^{-m_0 u^2/2}}{V} +
J^{(1)}_V(u) +J^{(2)}_V(u)+J^{(3)}_V(u) .
\end{eqnarray*}
We now have to choose $V$ in an appropriate way. The choice of $V$ will
be given by the bound on $J^{(1)}_V(u)$.
Assume in the following that $u$ satisfies the condition $(U1)$ given
by $\frac{u^2m_{22}}{2\lambda}\leq\frac{m_2}{4}$, and let us set
\[
V= \frac{3\sqrt{\lambda} m_2}{4m_{03}} .
\]
Then for all $v\in[0,V]$, $-\frac{v^2m_2}{2}+\frac{v^2}{2}\frac
{u^2m_{22}}{2\lambda} +\frac{|v|^3 m_{03}}{6\sqrt{\lambda}} \leq
-\frac{v^2m_2}{4}$, and thus
\[
J^{(1)}_V(u)\leq\frac{m_{03}}{3\sqrt{\lambda}} \int_{0}^{V} v e^{
-{v^2m_2}/{4}}  \,dv \leq\frac{2 m_{03}}{3 m_2 \sqrt{\lambda}}
.
\]
For the term $J^{(2)}_V(u)$, we notice that if $u$ satisfies the
condition $(U2)$ given by $\frac{|u|}{\sqrt{\lambda}}m_{12} \leq
\frac{m_2}{2}$, then for all $V>0$, we can bound $J^{(2)}_V(u)$ by
\[
J^{(2)}_V(u) \leq\frac{|u|}{\sqrt{\lambda}}m_{12}\int
_{0}^{V}e^{-({v^2}/{4})m_2}   \,dv \leq\frac{|u|}{\sqrt{\lambda
}}m_{12} \sqrt{\frac{\pi}{m_2}} .
\]
Finally, for the third term, we have that if $u$ satisfies the
condition $(U3)$ given by $\frac{|u| m_{30}}{3\sqrt{\lambda}}\leq
\frac{1}{2} m_0$, then $J^{(3)}_V(u)$ can be bounded, independently of $V$,~by
\[
J^{(3)}_V(u)\leq\sqrt{2\pi m_2}\frac{|u|^3 m_{30}}{6\sqrt{\lambda
}} e^{-({1}/{4}) u^2 m_0} \leq\sqrt{2\pi m_2}\frac{2 |u|
m_{30}}{3m_0\sqrt{\lambda}} e^{-1}
\]
because of the fact that for all $x\geq0$, then $xe^{-x}\leq e^{-1}$.

The final conclusion of all these computations is that if we set
$a_1=\min(\sqrt{\frac{m_2}{2m_{22}}},  \frac{m_2}{2m_{12}} , \frac
{3m_0}{2m_{30}})$, then for all $u$ and $\lambda>0$ we have
\[
|u|\leq a_1 \sqrt{\lambda} \quad\Longrightarrow\quad\bigl|\pi\widehat
{C_{Z_{\lambda}}}(u,[0,1]) - \sqrt{2\pi m_2} e^{-m_0 u^2/2}\bigr| \leq
\frac{a_2}{\sqrt{\lambda}} + \frac{a_3 |u|}{\sqrt{\lambda}} ,
\]
where $a_2 = \frac{24 m_{30}+ 2 m_{03}}{3m_2}$ and $a_3 = m_{12}\sqrt
{\frac{\pi}{m_2}} + \frac{2\sqrt{2\pi m_2} m_{30}e^{-1}}{3m_0}$.

\subsection{\texorpdfstring{Proof of Lemma \protect\ref{PathOfZeroslem}}{Proof of Lemma 3}}

The proof of this lemma relies upon the implicit function theorem. Let
us start with the proof of (i);
let $(\sigma_0,t_0)$ be a point such that $h(\sigma_0,t_0) =
0$. By Assumption (a), we have that $\frac{\partial h}{\partial t
}(\sigma_0,t_0)\neq0$. Then, thanks to the
implicit function theorem, there exist two open intervals $I=(\sigma
_0^{-} ,
\sigma_0^{+})$ and $J=(t_0^{-} ,t_0^{+})$ containing, respectively,
$\sigma_0$
and $t_0$, and a $C^1$ function $\eta\dvtx I\to J$ such that $\eta(\sigma_0)=t_0$
and $\forall(\sigma,t)\in I\times J$, $h(\sigma,t)=0 \Leftrightarrow
t=\eta(\sigma)$. Let us now denote $\eta=\Gamma_{t_0}$. We need to
prove that we can take $\sigma_0^{-}=0$ when $\Gamma_{t_0}$ remains
bounded. Assume we cannot; the maximal interval on which $\Gamma
_{t_0}$ is defined is $(\sigma_0^{-},\sigma_0^{+})$ with $\sigma
_0^{-}>0$. By assumption, there is an $M_0>0$ such that for all $\sigma
\in(\sigma_0^{-},\sigma_0^{+})$, then $|\Gamma_{t_0}(\sigma)|\leq
M_0$. We can thus find a subsequence $(\sigma_k)$ converging to
$\sigma_0^{-}$ as $k$ goes to infinity and a point $t_1\in[-M_0,M_0]$
such that $\Gamma_{t_0}(\sigma_k)$ goes to $t_1$ as $k$ goes to
infinity. By continuity of $h$, we have $h(\sigma_0^{-}, t_1)=0$. Now,
we also have $\frac{\partial h}{\partial t}(\sigma_0^{-},t_1)=0$.
Indeed, if it were $\neq0$, we could again apply the implicit function
theorem in the same way at the point $(\sigma_0^{-},t_1)$, and get a
contradiction with the maximality of $I=(\sigma_0^{-},\sigma_0^{+})$.
Then, by Assumption (b), we have $\frac{\partial h}{\partial\sigma
}(\sigma_0^{-},t_1)\neq0$.
We can again apply the implicit function
theorem, and we thus obtain that there exist two open intervals
$I_1=(\sigma_1^{-} , \sigma_1^{+})$ and $J_1=(t_1^{-} ,t_1^{+})$ containing,
respectively,\vadjust{\goodbreak} $\sigma_0^{-}$ and $t_1$, and a $C^1$ function $\xi
\dvtx J_1\to I_1$ such that
$\xi(t_1)=\sigma_0^{-}$ and $\forall(\sigma,t)\in I_1\times J_1$,
$h(\sigma,t)=0 \Leftrightarrow\sigma=\xi(t)$. Moreover, we can
compute the derivatives of
$\xi$ at $t_1$. We start from the implicit definition of $\xi$:
$h(\xi(t),t)=0$. By differentiation, we get $\xi'(t)\frac{\partial
h}{\partial\sigma} (\xi(t),t) + \frac{\partial h}{\partial t } (\xi(t),t)
=0$. Taking the value at $t=t_1$, we get $\xi'(t_1)=0$. We can again
differentiate, and find $\xi''(t)\frac{\partial h}{\partial\sigma
}(\xi(t),t) + \xi'(t)^2\frac{\partial^2 h}{\partial\sigma^2}(\xi
(t),t)+ 2\xi'(t)\frac{\partial^2 h}{\partial\sigma\,\partial t}(\xi
(t),t)+ \frac{\partial^2 h}{\partial t^2}(\xi(t),t)=0$. Taking again
the value at $t=t_1$, we get
\[
\xi''(t_1) = - \frac{1}{\xi(t_1)} = - \frac{1}{\sigma_0^{-}} < 0 .
\]
Thus it shows that $\xi$ has a strict local maximum at $t_1$; there
exist a
neighborhood $U_1$ of $\sigma_0^{-}=\xi(t_1)$ and a neighborhood
$V_1$ of $t_1$ such that for all points in $U_1\times V_1$, then
$h(\sigma,t)=0$ implies $\sigma=\xi(t)\leq\xi(t_1)=\sigma_0^-$,
which is in contradiction with the definition of $\Gamma_{t_0}$ on
$(\sigma_0^{-},\sigma_0^{+})$. This ends the proof of~(i), and also
of (iii).

For (ii), assume that $t_0$ and $\widetilde{t_0}$ are two points such
that $h(\sigma_0,t_0)=h(\sigma_0,\widetilde{t_0}) =0$ and such that
there exists $\sigma_1<\sigma_0$ such that $\Gamma_{t_0}(\sigma
_1)=\Gamma_{\widetilde{t_0}}(\sigma_1)=t_1$. Then, if $\frac
{\partial h}{\partial t}(\sigma_1,t_1)\neq0$, the implicit function
theorem implies that $\Gamma_{t_0}(\sigma)=\Gamma_{\widetilde
{t_0}}(\sigma)$ for all $\sigma\in[\sigma_1,\sigma_0]$ and in
particular $t_0=\widetilde{t_0}$. But now, if $\frac{\partial
h}{\partial t}(\sigma_1,t_1)=0$, then, as above, this implies that
$\frac{\partial h}{\partial\sigma}(\sigma_1,t_1)\neq0$ and using
again the implicit function theorem, this would be in contradiction
with the fact $\Gamma_{t_0}(\sigma)$ is defined for $\sigma\in
[\sigma_1,\sigma_0]$.
\end{appendix}

\section*{Acknowledgments}
The authors are also very grateful to the anonymous referees for their
careful reading
and their relevant remarks contributing to the improvement of the manuscript.


%


\printaddresses

\end{document}